\documentclass[11pt]{article}
\usepackage{amsmath}
\usepackage{amssymb}
\usepackage{txfonts}

\newtheorem{deff}{Definition}

\newtheorem{proposition}{Proposition}
\newtheorem{example}{Example}
\newtheorem{theorem}{Theorem}

\newcommand{\proof}{\vskip2mm \noindent {\bf Proof.}~}

\newcommand{\Aut}{\mathop{\mathrm{Aut}}}

\newcommand{\mto}{\mapsto}

\newcommand\beqa {\begin{eqnarray}}
\newcommand\eeqa {\end{eqnarray}}
\newcommand\bqa {\begin{eqnarray}}
\newcommand\eqa {\end{eqnarray}}
\newcommand{\beq}{\begin{eqnarray}}
\newcommand{\beqn}{\begin{eqnarray}\nonumber}
\newcommand{\eeq}{\end{eqnarray}}
\newcommand{\be}{\begin{array}}
\newcommand{\ee}{\end{array}}

\newcommand\bea {\begin{eqnarray}}
\newcommand\eea {\end{eqnarray}}
\newcommand\nn{\nonumber}

 \newcommand{\pr}{\partial}
 \newcommand{\pt}{\partial}

 \newcommand{\cE}{{\cal E}}

 \newcommand{\cV}{{\cal V}}
 
 \newcommand{\cP}{{\cal P}}

 \newcommand{\cM}{{\cal M}}

\newcommand{\R}{{\mathbb R}}
\newcommand{\Z}{{\mathbb Z}}

 \newcommand{\g}{{\mathfrak g}}

\newcommand{\md}{\mathrm{d}}

 \newcommand\bphi {\overline {\phi }}
 
\newcommand\bA {\overline {A }}
 \newcommand\bF {\overline {F} }
 \newcommand\bE {\overline {E} }

\newcommand{\rank}{\mathop{\mathrm{rank}}}
 \def\S{{\Sigma}}
\def\2{{\textstyle\frac{1}{2}}}
\def\ba{\begin{eqnarray}}
\def\ea{\end{eqnarray}}

   \def\CM{{\cal M}}

\def\CP{{\cal P}}

\DeclareMathOperator\Sim{\mbox{\Large$\sim$}}

\def\EF{{\,\!^E{}\!F_\phi}}
\def\Emd{\!\,{}^E\!\md}
\def\EL{\,\!{}^E\!L}
\def\ELi{\,\!{}^{E_i}\!L}
\def\ELo{\,\!{}^{E_1}\!L}
\def\ELt{\,\!{}^{E_2}\!L}
\def\EPhi{\,\!{}^E\!\Phi}
\def\Ephi{\,\!{}^E\!\phi}
\def\EcV{\,\!{}^E\!\cV} 
\def\EdPhi{\,\!{}^E\!\delta \Phi}
\def\EdePhi{\,\!{}^E\!\delta_\e \Phi}
\def\EdetPhi{\,\!{}^E\!\delta_{\e_2} \Phi}
\def\EdeoPhi{\,\!{}^E\!\delta_{\e_1} \Phi}
\def\EdzePhi{\,\!{}^E\!\delta^0_\e \Phi}

\newcommand{\phigra}{\phi^{\rm gra}}
\newcommand{\Phigra}{\Phi^{\rm gra}}
\newcommand{\Alt}{\mathop{\rm Alt}}

\def\o{\omega}
\def\O{\Omega}
\def\OE{\Omega_E^\cdot(M)}
\def\OEo{\Omega_{E1}^\cdot(M_1)}
\def\OEt{\Omega_{E2}^\cdot(M_2)}
\makeatletter
\def\bard{\protect\@bard}
\def\@bard{%
\relax
\bgroup
\def\@tempa{\hbox{\raise.73\ht0
\hbox to0pt{\kern.4\wd0\vrule width.7\wd0
height.1pt depth.1pt\hss}\box0}}%
\mathchoice{\setbox0\hbox{$\displaystyle\mathrm{d}$}\@tempa}%
{\setbox0\hbox{$\textstyle \mathrm{d}$}\@tempa}%
{\setbox0\hbox{$\scriptstyle \mathrm{d}$}\@tempa}%
{\setbox0\hbox{$\scriptscriptstyle \mathrm{d}$}\@tempa}%
\egroup
}
\makeatother

\makeatletter
\def\barb{\protect\@barb}
\def\@barb{%
\relax
\bgroup
\def\@tempa{\hbox{\raise.73\ht0
\hbox to0pt{\kern-.1\wd0\vrule width.7\wd0
height.1pt depth.0pt\hss}\box0}}%
\mathchoice{\setbox0\hbox{$\displaystyle\mathrm{b}$}\@tempa}%
{\setbox0\hbox{$\textstyle \mathrm{b}$}\@tempa}%
{\setbox0\hbox{$\scriptstyle \mathrm{b}$}\@tempa}%
{\setbox0\hbox{$\scriptscriptstyle \mathrm{b}$}\@tempa}%
\egroup
}
\makeatother

\begin{document}
  
  \def\sp{\mathfrak sp}
   \def\sll{\mathfrak sl}
   \def\P{{\mathbb P}}
   \def\H{\mathbb H}
   \def\a{\alpha}
   \def\b{\beta}
    \def\d{\delta}
    \def\t{\theta}
   \def\de{\delta}
   \def\la{\lambda}
   \def\e{\epsilon}
   \def\i{\imath}
   \def\ve{\varepsilon}


\begin{center}
{ \bf \Large Lie algebroid morphisms,
Poisson Sigma Models, \\ and off-shell
closed gauge symmetries}
\end{center}

\smallskip

\begin{center}\sl\large
 Martin Bojowald$^a$\footnote{e-mail
address: {\tt mabo@aei.mpg.de}}, Alexei
Kotov$^b$\footnote{e-mail adress: {\tt A.Kotov@tpi.uni-jena.de}},
Thomas Strobl$^b$\footnote{e-mail address: {\tt
Thomas.Strobl@tpi.uni-jena.de}}
\end{center}


   \begin{center}{\sl $^a$ Max-Planck-Institut f\"ur
   Gravitationsphysik, Albert-Einstein-Institut\\
   Am M\"uhlenberg 1, D-14476 Golm, Germany \\
  \vspace{0.5em}
$^b$ Institut f\"ur
Theoretische Physik, Universit\"at Jena, D--07743 Jena, Germany}
\end{center}

\vskip 10mm

{May 14, 2004}


\centerline{\large\bf Abstract}\vskip 5mm
Chern-Simons  gauge theories in 3 dimensions and the Poisson Sigma
 Model (PSM) in 2 dimensions are examples of the same theory, if their
 field equations are interpreted as morphisms of Lie algebroids and their
 symmetries (on-shell) as homotopies of such morphisms. We point out
 that the (off-shell) gauge symmetries of the PSM in the literature are
 not globally well-defined for non-parallelizable Poisson manifolds and
 propose a covariant definition of the off-shell gauge symmetries as left
 action of some finite-dimensional Lie algebroid.

 Our approach allows to avoid complications arising in the infinite
 dimensional super-geometry of the BV- and AKSZ-formalism. This preprint
 is a starting point in a series of papers meant to introduce
 Yang-Mills type gauge theories of Lie algebroids, which include and
 generalize the standard
 YM theory, gerbes, and the PSM.

\vskip 4.5mm

\section{Introduction}
Yang-Mills (YM) gauge theories are an important ingredient in our
present-day understanding of fundamental forces. 
On the mathematical side they are governed by a principal
fiber bundle $\pi \colon P \to \Sigma$, where $\Sigma$ is our
spacetime manifold. Here any fiber $\pi^{-1}(x)$, $x \in
\S$, is a $G$-torsor, where $G$ is the Lie structure  group of $P$. In
the simplest case when $P$ is a trivial bundle, $P \cong \Sigma \times
G$,  $\pi$ is the projection to the first factor, and the $G$-action
on $P$ is defined by right multiplication in the second
factor. For the standard model of elementary particle physics $G=SU(3)
\times SU(2) \times U(1)$, but also other, ``larger'' Lie groups come into
mind in the context of a further unification of fundamental
interactions. 

Gauge bosons correspond to connections in $P$, matter fields are
sections in associated fiber bundles (usually vector bundles),
and local gauge symmetries are the vertical automorphisms $\Aut_v(P)$ of $P$.
In the case of a trivial bundle the gauge bosons are just $\g$-valued
1-forms $A=A^I b_I$ on $\S$, where $\g$ is the Lie algebra of the gauge or
structure group $G$, $b_I$ is some basis in $\g$, and $I = 1, \ldots
\dim{\g}$. Sections of vector bundles then correspond to vector-valued
functions (or spinors) on $\S$  and the infinite dimensional group of
gauge transformations  $\Aut_v(P)$ just becomes isomorphic to
$\mathrm{Map}(\S,G)$. Thus infinitesimally  local gauge symmetries
are parametrized by $\e=\e^I b_I \in \mathrm{Map}(\S,\g)$ and one has
$\delta_\e A^I = \md \e^I + C^I_{JK} A^J \e^K$, where $C^I_{JK}$
denote the structure constants of the Lie algebra $\g$. 

All fundamental interactions fit into this framework except for
gravity. Even though it is possible to cast
general relativity in the language of a gauge theory of connections
\cite{AA}, the local gauge symmetries contain the diffeomorphisms of
$\S$. The Lie algebra of $\mathrm{Diff}(\S)$ consists of vector fields
on $\S$. This has to be contrasted to elements of $\mathrm{aut}_v(P)
\equiv
\mathrm{Lie} (\Aut_v(P))$, which always have a trivial projection to
$T\Sigma$. On the level of a Hamiltonian formulation of the theory
this usually leads to structure functions in the algebra of
constraints,
whereas for YM gauge theories
the algebra of constraints is 
governed by the structure
constants $C^I_{JK}$ (cf., e.g., \cite{HT} for details). Structure
functions of first class constraints are a typical feature of a
formulation of a theory with an open algebra of local symmetries,
where the commutator of infinitesimal
local or gauge symmetries closes only on-shell, i.e.~upon use of the field
equations. In YM theories, on the other hand, gauge symmetries always
form a closed algebra. In a
way, within YM theories many considerations
of local symmetries can be reduced to a finite dimensional group, the
structure group $G$, whereas for gravitational theories all of the 
infinite-dimensional group of local symmetries seems unavoidable. This
may be regarded as maybe one of the main obstacles in a successful
quantization of gravity along the lines of YM-gauge theories. 

It may be an important step to broaden the framework of YM-gauge
theories in such a way that also some gauge theories with an open
algebra of gauge transformations
fit into it, while still many considerations 
can be reduced to a purely finite-dimensional setting. In
\cite{Jenatalk} (cf.~also \cite{LAYM})
a particular program in this
direction has been proposed. Essentially,
the structural Lie group $G$ of a YM theory
is  replaced by (or generalized to) a so-called Lie groupoid; correspondingly,
the Lie algebra $\g$ generalizes to a so-called Lie algebroid
$E$.\footnote{Among others a Lie algebroid is a
vector bundle $E \to M$ carrying a Lie bracket for its sections; for
$M$ a point one obtains a Lie algebra, while $E=TM$ with the
Lie-Jacobi bracket for vector fields on $M$ is another
prominent example. We will recall the notion of a Lie algebroid in
the subsequent section; for further background material on Lie
algebroids and Lie groupoids we refer to the monograph
\cite{SilvaWeinstein} and references therein.} The present paper is
the first one in a series of papers devoted to this subject and aims
at providing part of the mathematical basis for the others. 

{}From some other perspective our goal is 
to provide a better understanding
and definition of ``non-linear gauge theories'', as they have been suggested
 already quite some time ago by van Nieuwenhuizen and collaborators, cf., e.g.,
\cite{nonlinear}. Heuristically, in such a theory one wants to replace
the structure constants $C^I_{JK}$ of a standard YM-theory by some
field-dependent quantities, which then generically will lead to a
theory with an open algebra of local symmetries, due to the
transformation of $C$s. In our approach, $C^I_{JK}$ will be the
structure functions of a Lie algebroid $E \to M$. $M$ then serves as a
target space for a Sigma Model so that the map $\S \to M$ locally
corresponds to a set of scalar fields $X^i(x)$ and the coefficients
$C^I_{JK}$ depend on these fields in general; from a physical perspective
these fields can be some kind of Higgs fields or, as shown in 
\cite{LAYM}, they can turn out to
be just some auxiliary fields that do not carry any propagating 
degrees
of freedom, but serve as moduli parameteres. 
In addition to
them locally one still has a set of 1-form gauge fields $A^I$.

In two spacetime dimensions, $\dim \S = 2$,  a prototype of such a
non-linear gauge theory is provided by the Poisson Sigma
Model (PSM) \cite{PSM1,Ikeda}. It is worth mentioning here
that in this particular spacetime dimension, essentially all
possible YM gauge theories and 2d  gravity theories
find a unifying formulation as particular PSMs (cf., e.g.,
\cite{PSM2,PartI}). In all of our work we want to use the PSM 
as a kind of main guiding example for developing a more general
theory. In  particular, in the present first paper we show how the
field equations and the gauge symmetries of this model are related to Lie
algebroids. Focusing on the case corresponding to a trivial principal
bundle, the field content of the PSM, locally described by a set of
couples $(X^i,A_i)_{i=1}^n$ of scalar and 1-form fields, respectively,
corresponds 
to vector bundle morphisms $\phi \colon T\Sigma \to T^*M$. Since $M$ is a
Poisson manifold, both the source and target vector bundle carry Lie
algebroid structures. The content of the field equations will then be
shown to be equivalent to requiring $\phi$ to respect the Lie
algebroid structures, i.e.~to be Lie algebroid morphisms. 

Whereas for Lie algebras it is very straightforward to define the
notion of a morphism, for Lie algebroids the situation is somewhat
more intricate. After setting the notation and collecting some
background material in section 2,
in section 3 several formulations of such a morphism will be
mentioned and related to one another.
Essentially one needs to dualize
the map $\phi$, requiring it to be an appropriate chain map.
However, in our final formulation, using the graph of $\phi$,
this can also be circumvented.

An important observation in this context is
that also YM-type gauge theories such as the Chern-Simons theory
fit into that framework: Flatness of a connection $A=A^I b_I$ in a trivial
principal fiber bundle is tantamount to the condition that the
corresponding map from $T\S \to \g$, $\xi \mapsto  A^I(\xi) b_I$, is a
Lie algebroid morphism. Correspondingly, in our investigations 
we will replace $T^*M$ of the
PSM by an arbitrary Lie algebroid $E_2$. In fact, for means of generality we
will also generalize $T\S$ to an arbitrary Lie algebroid $E_1$, although
the main example of physical interest
may still be provided by the tangent bundle of spacetime.

For the formulation of  $\phi\colon E_1 \to E_2$ in terms of the graph map
one uses the fact that the set  $E_1 \times E_2$ can be given
the structure of a Lie algebroid $E=E_1 \boxplus E_2$ itself (details
of this will be provided in section 2 below). It will then be
shown that $\phi$ is a morphism of Lie algebroids, iff
$\phigra \colon E_1 \to E$ is a morphism. By construction,
the base map of $\phigra$ is an embedding, permitting 
to work with $\phigra$-related sections instead of with
the dual map. 

In section 4 finally we turn to the issue of local gauge symmetries.
We first point out that the local infinitesimal 
symmetries usually used in the PSM are in general not well-defined
globally. They make sense only if the target Poisson manifold $M$
can be covered by a single chart, or if it carries some flat
connection, implicit but not transparent in the usual formulas
(Eqs.~(\ref{sym1}) and (\ref{sym2}) below). This is somewhat
remarkable in view of the already relatively large, and in part
also mathematical literature on the PSM; partially this may be related to
the fact that in many physical examples of the PSM
such as 2d YM- and/or 2d gravity models a flat target $M \cong
\R^n$ is used (cf., e.g., \cite{PartI,Habil2,Kummerreview}), which
moreover also underlies the Kontsevich formula \cite{Kontsevich},
resulting from the perturbative quantization of the PSM \cite{CF1}.

In section 4 we
present one possible way of curing this deficiency, simultaneously generalizing
the local symmetries also to the context of arbitrary Lie
algebroids. This is done in such a way that
for the particular case $E_2 := \g$ and  $E_1 := T\S$ 
one indeed reobtains the usual YM gauge transformations.
Moreover, also in the general case,
we will be able to trace back everything to purely finite
dimensional terms. Employing the picture with the graph,
$\phigra \colon E_1 \to E$, the infinitesimal
gauge symmetries (and also what
corresponds to infinitesimal diffeomorphisms of $\S$) result from
particular, structure preserving infinitesimal automorphisms of $E$,
acting from the left on $\phigra$ (or from the right in the
dualized picture $\Phigra \colon \Gamma(\Lambda^\cdot E^*) \to
\Gamma(\Lambda^\cdot E_1^*)$), and generated by particular
sections of $E$ via a Lie algebroid generalization of the Lie derivative.  
As a byproduct we find that the gauge symmetries formulated in this
way close even off-shell. But also if one needs to calculate e.g.~the
commutator of the original symmetries of the PSM for $M \cong \R^n$
the present approach provides a significant technical advance.

Although this approach may be related also to an infinite-dimensional
Lie algebroid $\cE$ of infinitesimal gauge transformations  \cite{LO},
the base manifold $\cM$ of which are maps $\S \to M$ (or, more
generally, maps from the base of $E_1$ to the base of $E_2$), one can
consistently---and with conceptual profit---truncate
$\Gamma(\cE)$ to the space of sections in the finite dimensional
algebroid $E$. For the PSM a likewise statement applies to its
AKSZ-formulation \cite{AKSZ,CFAKSZ}, which yields in a most transparent
way the BV-form of the PSM. 

As an alternative, one may also employ a connection in the target Lie
algebroid $E_2$ for providing another possible global definition of
the local gauge symmetries. While some elementary formulas in this
direction will be displayed at the end of section 4,
a more abstract analysis along the
lines of the present paper can be found in another, accompanying
paper \cite{AT1}.

Both definitions of gauge symmetries can be made to agree
for the PSM on $M\cong \R^n$, as well as certainly in the YM-case.
Also they always agree globally
upon use of the field equations, i.e.~on-shell. Already the standard
gauge symmetries of the PSM have a good global on-shell meaning, as an
infinitesimal homotopy of Lie algebroids. Correspondingly, a homotopy
of Lie algebroids defines an integrated version of the on-shell gauge
symmetries (section 4). Globally and
off-shell, however, the gauge symmetries defined via an $E$-Lie
derivative and those defined by means of a connection $\Gamma$ on $E_2$ are
different; in particular also the latter do not close off-shell, their
commutator containing contributions of $\Gamma$.

The formulation in the present paper as well as in \cite{AT1} is
put in such a form that a generalization to non-trivial fibrations is rather
straightforward. Essentially $E$, as a
manifold, is then not just a direct product $E_1 \times E_2$, but a particular
fiber bundle over the base of $E_1$. In order to not overload the
presentation, we found it useful to present this generalization in
another separate work \cite{AT2}. All three papers
together then are meant to provide, among others, 
a basic mathematical framework for the definition
of Lie algebroid Yang-Mills type gauge theories. 

Some examples for action functionals of this kind of gauge theories are
presented in \cite{LAYM}. They 
generalize e.g.~usual YM gauge theories in arbitrary dimensions of $\S$. It
turns out that one can use such theories to effectively glue together
YM theories with different structure groups, which can even differ in
their dimension, this theory being then governed by one Lagrangian;
due to the gauge symmetries, the
map from $\S \to M$, corresponding to the scalar fields $X^i(x)$,
carries only global degrees of freedom, and the representative maps
have the role of moduli for changing from one YM theory to another
one. But also topological action functionals can be constructed,
generalizing e.g.~the Chern-Simons gauge theory in three and the PSM in two
spacetime dimensions. Also one can extend the relation of the PSM to
2d gravity theories, 
to the definition of
topological gravity theories in arbitrary spacetime dimensions; for
this cf.~\cite{CMP1}. 
Still further work is necessary to see how far one can push this
approach and what kind of different theories can be constructed.
The present paper is meant to set part of the basis
for it. One may expect, however, that 
the resulting theories have a wide range of implications, in physics as
well as in mathematics.

\section{Preliminaries}
\label{sec:prelim}

In this section we mainly set the notation
and recall some background material needed later on. We
start with the Poisson Sigma Model (PSM) \cite{PSM1,Ikeda}, presenting
a slightly more abstract definition of its action functional $S$. $S$
is a functional of the vector bundle morphisms $\phi \colon T
\Sigma \to T^*M$, where $\Sigma$ is a two-dimensional manifold, called
the world sheet, and $M$ some Poisson manifold. We denote the Poisson
bivector by $\cP \in \Gamma(\Lambda^2 TM)$,  $\{ f,g \} = \langle \cP ,
\md f \wedge \md g \rangle$; in local coordinates $X^i$
on $M$, $\cP = \frac{1}{2} \cP^{ij}(X) \partial_i \wedge \partial_j
\Rightarrow \{X^i , X^j \} =  \cP^{ij}$, 
and
\begin{equation} [\cP, \cP ]_{\mathrm{Schouten}}
 \equiv 
 \cP^{ij}{}_{,s} \,\cP^{ks}  \,\partial_i \wedge \partial_j \wedge
 \partial_k =0,  \label{Jacobi} \end{equation}
as a manifestation of the Jacobi identity for the Poisson bracket. 

Any morphism $\phi \colon E_1 \to E_2$ between two
vector bundles $\pi_i \colon E_i \to M_i$, $i=1,2$, may be
expressed in different equivalent ways. One of them is by specifying
the induced base map $\phi_0 \colon M_1 \to M_2$ and, in addition, by
providing a section $A$ of the bundle $E^*_1 \otimes \phi_0^* E_2$. 
If $b_I$, $I=1, \ldots \rank(E_2)$, 
 denotes a local basis of $E_2$ and $\barb_I$ the corresponding
 induced basis in  the pullback
bundle $\phi_0^* E_2$, and if $E_1 = T\S$, then $A = A^I \otimes \barb_I$,
where $A^I \in \Omega^1(\Sigma) \equiv \Gamma(T^* \Sigma)$ (possibly
also defined locally on $\Sigma$ only, however).

Later on we will also need the graph $\phigra$
of the above map $\phi$ as well as
its trivial extension ${\,\!{}^E{}\phi}$,
\beq  \phigra \colon E_1 \to E := E_1 \boxplus E_2  \, , \quad e_1
\mapsto e_1 \boxplus \phi(e_1) \, , \label{phigra} \\
 {\,\!{}^E{}\phi} \colon E \to E \, , \quad {\,\!{}^E{}\phi} = \phigra \circ p_1
 \, .\label{phiext}
\eeq
Here $\pi \colon E\to M$, the exterior sum of $E_1$ and $E_2$,
is a  vector bundle over $M := M_1 \times
M_2$ defined as $\mathrm{pr}_1^* E_1 \oplus \mathrm{pr}_2^* E_2$,
where $\mathrm{pr}_i : M\to M_i$ is the projection to the $i-$th factor
of the Cartesian product, and
$p_1$ is the canonical projection bundle morphism $E \to E_1$ covering
$\mathrm{pr}_1 : M\to M_1$:
\beq
\be{ccc}
E & \stackrel{p_1}{\longrightarrow} & E_1 \\
\pi \downarrow & & \downarrow \pi_1 \\
M & \stackrel{\mathrm{pr}_1}{\longrightarrow} & M_1
\ee
\eeq 

Alternatively, the vector bundle morphism $\phi$ induces a map $\Phi
\colon \Gamma(\stackrel{p}{\otimes}\! E_2^*) \to
\Gamma(\stackrel{p}{\otimes}\! E_1^*)$. For $p=0$ it is given by the
pullback of functions, $C^\infty(M_2) \ni f \mapsto \phi_0^*f \in
C^\infty(M_1)$, while for $p=1$, $\Phi(u_2)$ for $u_2 \in
\Gamma(E^*_2)$ is defined by $\langle\Phi(u_2), s_1 \rangle|_x =
\langle u_2|_{\phi_0(x)}, \phi(s_1|_x) \rangle$ $\forall x \in M_1$
and $\forall s_1 \in \Gamma(E_1)$.  In the particular case of $E_1 =
T\S$ mentioned previously, with $b^I$ denoting the local basis in
$E_2^*$ dual to $b_I$, one has $\Phi(b^I) = A^I$.  The extension to
arbitrary $p$ is canonical now.  Mostly we will use only the
restriction of the above map $\Phi$ to the antisymmetric subspace
$\Gamma(\Lambda^p E_2^*)=:
\Omega^p_{E_2}(M_2)$ (the space of $E_2$-forms) only, which
we denote by the same letter. 

The above map $\Phi$ can be extended also to all $E_2$-tensors, and we
will denote this extension by
\beqa \label{phi!}
\Phi^! \colon \Gamma(\stackrel{p}{\otimes}\! E_2^*
\stackrel{q}{\otimes}\! E_2) \to
 \Gamma(\stackrel{p}{\otimes}\! E_1^* \stackrel{q}{\otimes}\!
 (\phi_0)^* E_2) \, , 
 \eeqa
where on the first factor $\Phi^!$ acts as $\Phi$ above and on $E_2$
it is defined as $\Gamma(E_2) \ni s_2 \mapsto s_2 \circ \phi_0$,
viewed as a section of the pullback bundle $(\phi_0)^*E_2$. With this
map the above section $A \in E^*_1 \otimes \phi_0^* E_2$ is nothing
but the image of the canonical identity section $\delta \in E_2^*
\otimes E_2$, $A = \Phi^!(\delta)$ (in local terms $\delta = b^I
\otimes b_I$ and $\Phi^!(b_I) \equiv \barb_I$).

In the particular case $E_2 =T^*M$ and $E_1 =T\S$ (and only in this case!)
the map $\Phi$ can be extended to all $E_2$-tensors also in another
way, which we denote by 
\beq\label{Phistar}\Phi^* \colon
\Gamma(\stackrel{p}{\otimes} TM \stackrel{q}{\otimes} T^*M ) \to
\Gamma(\stackrel{p+q}{\otimes} T^*\S)\,. \eeq
Here 1-forms on $M$, corresponding to $p=0$, $q=1$, are mapped by the
pullback $\phi_0^*$ to 1-forms on $\S$---and, as before, this map is
extended canonically to all possible choices for $p$ and $q$. 

Such as $\phi$ permits the dual formulation in terms of $\Phi
\colon \Omega^p_{E_2}(M_2) \to \Omega^p_{E_1}(M_1)$, also the maps
(\ref{phigra}) and (\ref{phiext}) induce reverse maps:
\beq \Phigra \colon \Omega^p_E(M) \to  \Omega^p_{E_1}(M_1) \; ,
\quad 
{\,\!{}^E{}\Phi} \colon \Omega^p_E(M)\to \Omega^p_E(M)  \, . \label{Phimaps}
\eeq 
Note that due to the isomorphism
$\bigoplus\limits_{p+q=k} \Omega^p_{E_1}(M_1) \otimes
\Omega^q_{E_2}(M_2) \cong 
\Omega_E^{k}(M) $, where multiplication is defined according to
$(\o_1 \otimes \o_2) \wedge (\o_1' \otimes \o_2') = (-1)^{qp'}
(\o_1  \wedge
\o_1') \otimes (\o_2 \wedge \o_2' )$, 
 there is a natural bigrading for $E$-forms;
if we want to stress the bigrading, we write $\Omega_E^{p,q}(M)$,
while $k$ in 
$\Omega_E^{k}(M)$ 
denotes the total degree, which is the sum of the
two individual degrees on $E_1$ and $E_2$. The above maps
(\ref{Phimaps}) are related to $\Phi$ in the following way: 
\beq \label{Phigra}
\Phigra \colon \Omega_E^{p,q}(M) \cong
\Omega^p_{E_1}(M_1) \otimes \Omega^q_{E_2}(M_2)
\stackrel{\mathrm{id} \otimes \Phi}{\longrightarrow}
\Omega^p_{E_1}(M_1) \otimes \Omega^q_{E_1}(M_1)
\stackrel{\wedge}{\rightarrow} \Omega^{p+q}_{E_1}(M_1)
  \, ,
\\ \label{PhiE}
{\,\!{}^E{}\Phi} \colon \Omega_E^{p,q}(M)
\stackrel{\Phigra}{\longrightarrow} \Omega^{p+q}_{E_1}(M_1)
\stackrel{P_1}{\longrightarrow}  \Omega^{p+q}_{E_1}(M_1)
\otimes \Omega^0_{E_2}(M_2) \cong  \Omega_E^{p+q,0}(M) \, , 
\eeq
where $P_1 \colon \O^p_{E_1}(M_1) \ni \o_1 \mapsto
\omega_1 \otimes 1 \in \O^p_{E_1}(M_1) \otimes \O^0_{E_2}(M_2) \subset
 \O^p_{E}(M)$ is the map induced by
the bundle morphism $p_1 \colon E \to
E_1$.
So, ${\,\!{}^E{}\Phi}$ preserves only the total
degree, but not the bigrading. 

By definition, $p_1 \circ \phigra = \mathrm{id}_{E_1}$, which
translates into the dual relation $\Phigra \circ P_1 =
\mathrm{id}$. For  ${\,\!{}^E{}\Phi}=P_1 \circ \Phigra$ we then
obtain ${\,\!{}^E{}\Phi} \circ P_1 = P_1$, implying that 
${\,\!{}^E{}\Phi}$ is a projector to the image of $P_1$
(i.e.~on $\Omega_E^\cdot(M)$ one has
${\,\!{}^E{}\Phi}^2=\EPhi$ and $\mathrm{im}\EPhi =
\mathrm{im}P_1$).


\vskip7mm

Using the map $\Phi^*$, we may now give a concise global definition of
the action functional of the PSM:\footnote{In section
\ref{sec:Generalized} this formula is rewritten in two further
quite similar fashions, cf.~Eq.~(\ref{actiondelta2}) below, which
will be explained only there to not overload
the presentation here.}
\beq S[\phi] = \int_\S \Alt \, \Phi^* (\delta+\cP) \, , \label{actiondelta}
\eeq
where $\delta$ is the canonical identity section
in $TM \otimes T^*M$ and $\Alt$ denotes the antisymmetrization. In
local coordinates $X^i$ on $M$ and with the induced local basis
$\partial_i \sim b_I$
and $\md X^i \sim b^I$ in $TM$ and $T^*M=E_2$, respectively, one has $\delta =
\partial_i \otimes \md X^i$ and  $ \Alt \,
\Phi (\delta) = \Alt \, (A_i \otimes \md X^i) = A_i \wedge \md
X^i$ (where $A_i \sim A^I$, as introduced above, and $X^i=X^i(x)$
denotes the scalar field corresponding to the map $\phi_0 \colon \S
\to M$, just expressed in local coordinates). $\cP$, on the other hand, 
is the Poisson tensor on $M$, and for the
second term simply $ \Alt \, \Phi^* (\cP) = \Phi
(\cP) = \2 \cP^{ij} A_i \wedge A_j$. Thus in the more familiar and for
practical purposes most useful local description, $S$ takes the form
\beq S =S[\phi_0,A]= \int\limits_\Sigma A_i \wedge \md X^i +
\frac{1}{2} \cP^{ij} A_i \wedge
A_j \, .
\label{action} \eeq
For completeness we also mention another possible covariant presentation
of the action functional: For this purpose we first rewrite $\cP$ as 
$\frac{1}{2}\langle \CP , \delta \wedge \delta \rangle$, then
the second term in (\ref{actiondelta}), which may be also written as
 $\Phi^!(\cP)$, becomes $\frac{1}{2}
\langle \Phi^!(\CP) , A \wedge A \rangle$ with $A \wedge A \in
\Omega^2(\S,\Lambda^2 \phi_0^* T^*M)$.  Moreover, $(\phi_0)_*
\colon T\S \to TM$ is a vector bundle morphism covering $\phi_0$.
Thus, according to the above discussion, it induces a section of
$T^*\S \otimes \phi_0^* TM$, which we
denote suggestively by  $\md \phi_0$. It clearly can be contracted
with $A \in \Gamma(T^*\S \otimes  \phi_0^* T^*M)$. In this way we obtain
\beq S =
\int_\S \langle A \stackrel{\wedge}{,} \md \phi_0 \rangle + \frac{1}{2}
\langle \CP \circ \phi_0 , A \wedge A \rangle \, .
\label{actionprime}
\eeq

Concerning the field equations and the symmetries of the PSM
action functional, we let it suffice here to just recall the local basis
expressions---anyway, much to follow will be
devoted to a more abstract and covariant formulation of precisely
these two issues. 

The field equations of the action functional (\ref{action}) are
\beq \frac{\delta S}{\delta A_i} &\equiv& \md X^i + \cP^{ij}(X) A_j = 0
\label{eqs1} \\   \frac{\delta S}{\delta X^i} &\equiv& \md A_i + \frac{1}{2}
\cP^{kl}{},_{i}(X) \, A_k \wedge A_l = 0 \label{eqs2} \eeq
The gauge symmetries are generated by
\ba
\de_\e X^i &=& \cP^{ji}\e_j  \label{sym1} \\
\de_\e A_i &=& \md\e_i +\cP^{jk}{},_i A_j \e_k \label{sym2}
\ea 
where $\epsilon \equiv \epsilon_i \md X^i \in \Gamma(\phi_0^* T^*M)$
may be chosen arbitrarily. The obvious Diff($\S$) invariance of 
the action functional $S$, e.g., can be generated by means of (\ref{sym1}) and
(\ref{sym2}) with the choice $\e_i = \langle v,A_i \rangle$ with $v \in
\Gamma(T\Sigma)$ being the infinitesimal generator of a
diffeomorphism in the above group.  For further remarks in the context
of symmetries, somewhat complementary to what will follow in the
present paper, we also refer to the introductory section 2.1 of
\cite{PSMRiemann}.

\vskip7mm

We now recall the definition of a Lie algebroid. First of all,
$E=T^*M$, $M$ Poisson, is a particular example, and many things become
more transparent when they are formulated in this somewhat more
general context and language. Moreover, although the action functional
$S$, as introduced above, is quite particular to morphisms from only
$T\Sigma \to T^*M$, where $\S$ is two-dimensional and $M$ Poisson, the
field equations and symmetries generalize easily to arbitrary Lie
algebroid morphisms $\phi \colon E_1 \to E_2$. Moreover, we believe
that the corresponding considerations are of interest in this more
general context as well. Finally, we remark that it is even possible
to construct action functionals for this more general setting, too,
but this is not subject of the present paper.

A Lie algebroid over a base manifold $M$ is a vector bundle $E$ with a
Lie algebra structure $[\cdot,\cdot]$ on the space of sections $\Gamma (E)$
together with a bundle
map $\rho \colon E \to TM$, called the {\sl anchor}, which by
definition governs the following Leibniz rule: for any
$s, s'\in \Gamma (E)$, $f\in C^\infty (M)$,
\beq\label{Leibniz}
[s ,fs']=f[s , s'] +\rho_{s}(f)s' \, ,
\eeq
where $\rho_\cdot$ denotes the induced map of sections from
$\Gamma(E)$ to $\Gamma(TM)$. It is not difficult to see that 
$\rho_\cdot$ provides a representation of  $(\Gamma (E), [,])$ 
in the Lie
algebra of vector fields, i.e.~that
$[\rho_{s},\rho_{s'}]=\rho_{[s ,s']}$.
We briefly recall the list of standard examples of Lie algebroids: Lie
algebras, $M$ being a point, or bundles of Lie algebras, for $\rho
\equiv 0$. The tangent bundle, $E=TM$, $\rho = {\mathrm id}$. And,
finally, $E=T^*M$, $M$ Poisson, where $\rho = \cP^\sharp$,
$\rho(\alpha_i \md X^i) = \alpha_i \cP^{ij}
\partial_j$, and the bracket $[\md f, \md g] := \md \{f,g\}$
between exact 1-forms is extended to all 1-forms by means of
(\ref{Leibniz}).

There is also an equivalent definition of a Lie algebroid
$(E,M,\rho,[\cdot ,\cdot])$
as the differential graded algebra $(\Gamma (\Lambda^\cdot E^*), \wedge, \Emd)$,
where $\Emd$ is defined by 
($\omega \in \Gamma (\Lambda^p E^*)$, $s_i \in
\Gamma(E)$) 
\ba\label{Cartan}
\Emd \omega (s_1 ,...,s_{p+1})=\sum\limits_i (-1)^{i+1}
\rho (s_i)\o(..., \hat{s}_i,...)+
\sum\limits_{i<j} (-1)^{i+j}
\o([s_i ,s_j ], ...,\hat{s}_i,..., \hat{s}_j,...) ,
\ea
which is a generalization of the Cartan formula
 for the exterior derivative in the standard Lie algebroid $TM$.

An anchor map of a Lie algebroid $E$ provides a representation of
$\Gamma(E)$ in $C^\infty(M)$. 
One can lift this action to a
representation in  $\Gamma
(\Lambda^\cdot E^*)$: 
Taking any section $s$ of $E$, we associate a Lie
derivative ({\sl E-Lie derivative}) ${}^E\!L_s$ along $s$ by
generalization of 
Cartan's magic formula 
\beq\label{Lie_derivative}
{}^E\!L_s = [\Emd, \i_s ]=\Emd \i_s + \i_s \Emd , 
\eeq
where $\i_s$ denotes contraction with $s$ and $\Emd$ is defined in
(\ref{Cartan}) above. It is now straightforward to prove that indeed
one has a representation, i.e.~that $ [\,\!{}^E\!L_s , \EL_{s'}]=\EL_{[s
,s']}$ holds true.\footnote{This is done most easily by noting that 
the operator $[ \i_s ,[\i_{s'},\Emd]]$ on $\OE$ is
$C^\infty(M)$-linear and agrees with $\i_{[s,s']}$, cf.~also \cite{Kos,KS}.}
 (In general, for operators $\cV_1$,  $\cV_2$ 
of some fixed degree in a graded vector space, we define the
graded commutator bracket according to $[\cV_1,  \cV_2] := \cV_1 \circ \cV_2
- (-1)^{\deg \cV_1 \deg \cV_2} \cV_2 \circ \cV_1$. In the above,
$\Emd$, $i_s$, and $\,\!{}^E\!L_s$ are of degree $+1$, $-1$, and 0,
respectively.)



For later use we will need some of the above formulas in more explicit
form: Let $(U,\{X^i\})$ be a local coordinate chart, 
$b_I$ be a frame of $E_U$ over $U$, and  $b^I$ its dual frame in
$E^*_U$.
Then with 
$\rho (b_I) \equiv \rho_I =: \rho^{i}_I\pt_i$ and $[b_I,b_J] =:  C_{IJ}^K b_K$
 one finds 
\beq\label{can_target}
\Emd X^i = b^I \rho^{i}_I (X), \hspace{4mm}
\Emd b^I = -\frac{1}{2} C_{JK}^I (X) b^J \wedge b^K ,
\eeq
as well as 
\beq\label{Lie_2}
\EL_s X^i = s^I \rho^{i}_I , \hspace{4mm}
\EL_s b^I = \rho_J(s^I )b^J + C_{JK}^I (X) b^J s^K .
  \eeq
In the Poisson case, $b_I \sim \md X^i$, $b^I \sim \partial_i$,
$\rho^{j}_I \sim \CP^{ij}$,
 and $C_{JK}^I \sim \CP^{jk}{}_{,i}$.

\vskip7mm

Some words about conventions 
may be in place: If there are two Lie
algebroids involved, $E_i \to M_i$, $i=1,2$, such as already above in
the context of a bundle map $\phi \colon E_1 \to E_2$, we will mostly
mark objects of the respective algebroid with the corresponding
index. E.g.~$s_2, s_2' \in \Gamma(E_2)$ for sections of the target
bundle. Similarly, for the respective Lie algebroid exterior
derivatives, we will use the abbreviations ${}^{E_i}\md =: \md_i$.
However, to simplify notation we will make exceptions from the above
rule for what concerns e.g.~local coordinates and frames: $x^\mu,
b_\alpha$ denote coordinates and frame in the source $M_1$ and $E_1$,
respectively, while $X^i$ and $b_I$ do so for the
target. Correspondingly, then $C_{IJ}^K$ ($C_{\a\b}^\gamma$) denote
structure functions in $E_2$ ($E_1$), and likewise for connection
coefficients etc. Depending on the context, furthermore, $X^i$ may
just denote coordinates on $M_2$ or, as e.g.~already in (\ref{eqs1})
above, the collection of functions on (parts of) $M_1$ corresponding
to the base map $\phi_0 \colon M_1 \to M_2$; otherwise we would have
to write $\Phi(X^i) \equiv (\phi_0)^* X^i$, in the previously
introduced notation, where, moreover, $\Phi$ and $(\phi_0)^*$ are the
canonical restrictions of the respective maps to functions defined on
the neighborhood $U\subset M_2$ on which the coordinates $X^i$ are
defined. Likewise $\md X^i$ may denote a basis of local 1-forms in
$T^*M$ or its pullback, which more carefully we would have to write as
$(\phi_0)^*\md X^i\equiv \Phi^*\md X^i$. On the other hand, for the
induced basis in $(\phi_0)^*T^*M$ for clarity we use $\bard X^i :=
\Phi^! (\md X^i) \equiv \md X^i \circ \phi_0$. In generalization of
the 1-form fields $A_i$ of the PSM, we have the (locally defined) set
of $E_1$-1-forms $A^I \equiv A_{\alpha}^I \otimes b^\a = \Phi(b^I)$;
they combine into (the globally defined) $A = \Phi^!(b^I \otimes b_I)=
A^I \otimes \barb_I \in \Gamma(E^*_1 \otimes \phi_0^* E_2)$, which in
the PSM case becomes $A=A_i \otimes \bard X^i$.  



Finally we mention that if $E_i \to M_i$, $i =1,2$, are two Lie
algebroids, then also $E \to M$, where $E \equiv E_1 \boxplus E_2$  and
$ M \equiv M_1 \times M_2$ as introduced above, can be endowed
canonically with a  Lie algebroid structure (generalizing the direct sum
of two Lie algebras). For this purpose we use the
isomorphism $\Omega_E^\cdot (M)\cong
\Omega_{E_1}^\cdot (M_1) \otimes \Omega_{E_2}^\cdot (M_2)$, and define
$^E\md := {^E{}\md_1} + {^E{}\md_2} $, where
${^E{}\md_1} = \md_1 \otimes \mathrm{id}$ and, similarly,
${^E{}\md_2} =(-1)^{\ve_1}\mathrm{id}\otimes \md_2$,
with $\ve_1$ being the grading operator acting as multiplication
by $p$ on $\Omega_{E_1}^p (M_1)$. By construction, 
$({^E{}\md_i})^2=0$, and, due to the grading operator  $\ve_1$, also
${^E{}\md_1}$ and ${^E{}\md_2}$ anticommute, so that indeed
$(^E\md)^2=0$.


\section{Morphisms and field equations}
\label{sec:morph}
Assume that $E_1\to M_1$ and
$E_2\to M_2$ are Lie algebroids with the anchors $\rho_1$ and $\rho_2$,
respectively and that  $\phi\colon E_1 \to E_2$ is a vector bundle
morphism. For the particular case $E_1=T \Sigma$ and $E_2 = T^*M$, $M$
Poisson, $\phi$ reproduces the content of the fields in the PSM; it is
worthwhile, however, to discuss the more general situation
$\phi\colon E_1 \to E_2$ 
(cf.~also \cite{CMP1,LAYM} for further motivation for this perspective).
In the beginning of the
present section we address the question, under what
conditions we may call $\phi$  a morphism of Lie algebroids, as
well as how, in the particular case of the PSM, this is related to its
field equations. On our way we will prove also some helpful reformulations
of the notion of Lie algebroid morphisms in terms of the maps
introduced in the previous section. 

For $M_1 =M_2=\{pt\}$ the above Lie algebroids simply become Lie algebras.
By definition, $\phi\colon \g_1\to \g_2$ is a morphism of Lie algebras
iff  
$[\phi (s_1 ),\phi (s_1')]-\phi([s_1 ,s_1 '])=0$ $\;\forall 
s_1 , s_1 ' \in \g_1$. But, in general a vector bundle
morphism $\phi\colon E_1 \to E_2$ does not
induce a map of sections of those bundles (except if, say,
 the induced base map $\phi_0\colon M_1 \to
M_2$ is a diffeomorphism). Instead, as with vector fields and the
tangent map $\varphi_\ast$ of a map $\varphi \colon M_1 \to M_2$
(corresponding to the example of standard Lie algebroids $E_i =TM_i$
with $\phi =\varphi_\ast$),
 one may speak of relation of sections only: Sections $s_i \in
\Gamma (E_i)$ are called {\em $\phi$-related}, $s_1 \stackrel{\phi}{\sim}
s_2$, iff $\phi \circ s_1 = s_2 \circ \phi_0$.
Following  \cite{HigMc} we also say that $s_1\in \Gamma (E_1)$ is
$\phi-$projectable if it is $\phi-$related to some $s_2\in \Gamma (E_2)$.
The most straightforward attempt to generalize the morphism of Lie
algebras would then be
\begin{deff}
Let $\phi$ be a vector bundle morphism $\phi\colon E_1 \to E_2$ \label{weak} 
between two Lie algebroids \\ $(E_i,M_i,\rho_i,[\cdot , \cdot]_i)$,
$i=1,2$. We say that \emph{$E_1$ and $E_2$ are $\phi -$related},
$E_1 \stackrel{\phi}{\sim} E_2 $, iff
\begin{eqnarray} 
\rho_2 \circ \phi &=& (\phi_0)_\ast \circ \rho_1  \label{commute} \\  
s_1 \stackrel{\phi}{\sim}
s_2,  s'_1 \stackrel{\phi}{\sim}
s'_2 &\Rightarrow&  [s_1,s'_1]_1 \stackrel{\phi}{\sim} \label{relate}
[s_2,s'_2]_2 \qquad \forall s_i, s'_i \in \Gamma(E_i) 
\end{eqnarray}
where $(\phi_0)_\ast \colon TM_1 \to TM_2$ denotes the push forward of
tangent vectors induced by $\phi_0$. 
\end{deff} 
In general, however, $\phi$-relation of Lie algebroids is too weak a notion to
deserve being called also a morphism of Lie algebroids. We
thus take recourse to a dual perspective, using the map $\Phi$
introduced in the previous section (in the example of standard Lie
algebroids $E_i = TM_i$ and $\phi = \varphi_\ast$, the map $\Phi$ is just
 the pull back of differential forms): 
\begin{deff}\label{strong} 
A vector bundle morphism $\phi\colon E_1 \to E_2$ 
between two Lie algebroids \newline $(E_i,M_i,\rho_i,[\cdot ,
\cdot]_i) \simeq (\Gamma (\Lambda E_i^*), \wedge, \md_i)$,
$i=1,2$, is a \emph{morphism of Lie algebroids}, iff
the induced map $\Phi \colon \Gamma (\Lambda E_2^*) \to
\Gamma (\Lambda E_1^*)$
is a chain map:
\beq \label{LA_morph} \md_1 \Phi -\Phi\md_2 =0    \,. 
\eeq
\end{deff}
In other words, $\phi$ is a morphism iff
\beq  F_{\! \phi} \colon
\Omega^\cdot_{E_2}(M_2) \to \Omega^{\cdot + 1}_{E_1}(M_1) \; , \quad 
 F_{\! \phi} := \md_1 \Phi -\Phi\md_2  \label{Fphi}
\eeq
vanishes. 
Before continuation, we show that Definition \ref{strong} indeed 
serves the purpose of giving a mathematical meaning to the
field equations of the PSM:
\begin{proposition}\label{M2}
A bundle map $\phi$ between $T\Sigma$ and $T^*M$ is a solution of the PSM
equations (\ref{eqs1}, \ref{eqs2}), iff $\Phi$ is a morphism of
Lie algebroids.  
\end{proposition}
{\bf Proof.}
Let us choose a local chart $U\subset M$ supplied with coordinate
functions $\{X^i\}$, inducing the local frame $\pt_i$ of $TU$.
Applying $\md\Phi -\Phi\pt$ to $X^i$ and
$\pt_i$, we immediately obtain the first and the second field equations,
(\ref{eqs1}) and (\ref{eqs2}), respectively. Here $\md$ is the usual de Rham
operator on $\Sigma$ and $\pt$ is the Lichnerowicz-Poisson
 differential acting on
$\Gamma (\Lambda^\cdot TM)$, which is a particular case of the canonical
Lie algebroid
differential on $T^*M$ determined by the Poisson structure
$\cP$.  Since both the conditions (\ref{eqs1}), (\ref{eqs2}) and
(\ref{LA_morph}) are local, this completes the proof.
$\blacksquare$

In \cite{HigMc}, instead of the above, one finds the following
definition: 
\begin{deff} \label{HigMc} Let $E_1 ,E_2$ be Lie algebroids on bases
$M_1 ,M_2$ with anchors $\rho_1 ,\rho_2$. Then \emph{a morphism of Lie
algebroids} $E_1 \to E_2$ is a vector bundle morphism $\phi\colon
E_1\to E_2$, $\phi_0\colon M_1\to M_2$ such that equation
(\ref{commute}) holds
and such that for arbitrary
$s_1 ,s_1'\in \Gamma (E_1)$ with $\phi-$decomposition
\beq \label{decomp}
\phi \circ s_1 =\sum a_i (\eta_i \circ \phi_0 ), \hspace{4mm}
\phi \circ s_1' =\sum a_i' (\eta_i' \circ \phi_0 )
\eeq
we have
\beq \label{Higg}
\phi \circ [s_1 ,s_1' ] =\sum a_i a_j' ([\eta_i ,\eta_j']\circ \phi_0)+
\sum \rho_1 (s_1) ( a_j') (\eta_j' \circ \phi_0) - \sum
\rho_1 (s_1' )( a_j) (\eta_j \circ \phi_0) \, . 
\eeq
\end{deff}
Here $\{\eta_i\}, \{\eta_i'\}$ are sections of $E_2$ and $a_i ,a_i'$ 
functions over $M_1$. Let us mention that \emph{any} section $s \in
\Gamma(E_1)$ has \emph{some}  $\phi-$decomposition (e.g.~choose for
$\{\eta_i\}$ a (possibly overcomplete) basis of sections in
$E_2$---the definition then may be shown to be also independent of this
choice of basis).

\begin{proposition}
Definitions \ref{strong} and \ref{HigMc} are equivalent.
\end{proposition}
{\bf Proof.} As seen by a simple straightforward calculation,
application of (\ref{LA_morph}) to functions yields a dual formulation
of (\ref{commute}) (just contract the former equation with 
sections of $E_1$).

It remains to show equivalence of the second defining property in
Definition \ref{HigMc} to the application of (\ref{LA_morph}) to
sections of $E_2^*$. In other words we need to prove that for
any $u\in\Gamma (E_2^*)$ and $s_1
,s_1'\in\Gamma (E_1)$ with decompositions (\ref{decomp}) 
one has
\ba\label{mm1}
\langle (\md_1\Phi -\Phi\md_2 )u ,s_1\wedge s_1' \rangle =
\langle u\circ\phi_0 , \sum \rho_1 (s_1) ( a_j') (\eta_j'\circ\phi_0) -
\\ \nonumber -\sum
\rho_1 (s_1' )( a_j) (\eta_j \circ\phi_0 )-  
\phi\circ [s_1,s_1' ] +
\sum a_i a_j'([\eta_i ,\eta_j']\circ\phi_0 )\rangle \, .
\ea
In fact, using (\ref{Cartan}), we obtain
\ba
\langle \md_1\Phi (u) ,s_1\wedge s_1' \rangle =
\rho_1 (s_1) \langle \Phi u ,s_1' \rangle - \rho_1 (s_1') \langle
\Phi u, s_1 \rangle - \langle \Phi u, [s_1 ,s_1']\rangle = \\ \nonumber =
\rho_1 (s_1) ( \sum a_j' \phi_0^* \langle u, \eta_j'\rangle ) - 
\rho_1 (s_1' )( \sum a_j \phi_0^* \langle u, \eta_j \rangle )- 
\langle u \circ \phi_0 , \phi \circ [s_1,s_1' ]\rangle \, . 
\ea
The Leibniz rule for the anchor map action of $s_1 ,s_1'$ gives
\ba\label{mm2}
\langle \md_1\Phi (u) ,s_1\wedge s_1' \rangle =
\langle u\circ\phi_0 , \sum \rho_1 (s_1) ( a_j') (\eta_j'\circ\phi_0) 
 -\sum
\rho_1 (s_1' )( a_j) (\eta_j \circ\phi_0 )-
\phi\circ [s_1,s_1' ] \rangle +\\ \nonumber
+ \sum a_j' \rho_1 (s_1 )\phi^*_0 \langle u, \eta_j'\rangle -
\sum a_j' \rho_1 (s_1' ) \phi^*_0 \langle u,\eta_j\rangle\, .
\ea

On the other hand,
\ba\label{mm3}
\langle \Phi\md_2 u ,s_1\wedge s_1'\rangle =\sum a_i a_j'
\phi_0^* \langle \md_2 u, \eta_i\wedge \eta_j'\rangle =
\sum a_j a_j' \phi_0^* \left( \rho_2 (\eta_i )\langle u,\eta_j'\rangle
-\rho_2 (\eta_j' )\langle u,\eta_i\rangle
\right) - \\ \nonumber -
\langle u\circ\phi_0 , \sum a_i a_j' ([\eta_i ,\eta_j' ]\circ\phi_0)
\rangle
\ea
Eq.~(\ref{commute}) implies that
$\forall h\in C^\infty (M_2), s\in\Gamma (E_1), x\in M_1$, one has
$\rho_1 (s)_{|x}\phi^*_0 h = \rho_2 (\phi\circ s_x) h$.
Hence, taking into account the $\phi -$decompositions of $s_i ,s_j'$, we get
\beq
\rho_1 (s_1 )\phi^*_0 \langle u,\eta_j'\rangle =
\sum a_i \phi_0^* \rho_2 (\eta_i )\langle u,\eta_j'\rangle,
\eeq
and a likewise formula with primed and unprimed quantities
exchanged. Thus all additional contributions in the difference
$\langle \md_1\Phi u, s_1\wedge s_2\rangle -\langle\Phi\md_2 u ,
s_1\wedge s_1' \rangle$ vanish, i.e.~the last two terms in (\ref{mm2})
cancel against the first two terms in (\ref{mm3}).
$\blacksquare$

{}From  Definition \ref{HigMc} it is also obvious that
for Lie algebras, corresponding to $M_1 =M_2=\{pt\}$, the chain
property (\ref{LA_morph}) is equivalent to $\phi$ being a 
morphism in the usual sense. Also, from this version we see
that  if $\phi\colon E_1\to E_2$
is a morphism of Lie algebroids, then $E_1$ and $E_2$ are $\phi-$related:
Indeed, the condition on the left-hand side of (\ref{relate}) implies
a $\phi-$decomposition (\ref{decomp}) with only one term, $a=1$ and
$\eta = s_2$ (and likewise for the primed quantity), in which case
Eq.~(\ref{Higg}) just reduces to the right-hand side of (\ref{relate}).

However, in general the converse conclusion is not true as illustrated
e.g.~by the following example in the context of the PSM (cf.~our
discussion above and in particular Proposition \ref{M2}):
\begin{example}\label{counterexample} Let
$x^1,x^2$ be coordinates on the world-sheet $\S:=\R^2$ and
let $M:=\R^4$ be  a target manifold  supplied with a zero Poisson tensor.
Assume that $\phi$ is specified by the following choice of fields:
$A:=A_i \otimes \bard X^i$ with 
$A_1 := \md x^1$, $A_2 := x^2
\md x^1$, $A_3 := \md x^2$, $A_4 = x^2 \md x^2$ and  $\phi_0$, in accordance
with  the first morphism property (\ref{commute}), which is equivalent to
the first set of field equations (\ref{eqs1}), is chosen to map
to a single point in $\R^4$. This provides a
$\phi-$relation of $T\R^2$ and $T^*\R^4$, because there is not even a
single vector field $\xi$ on $\R^2$ that---for this choice of
$\phi$---is $\phi$-related to \emph{any} section of $\,\Gamma (T^*\R^4)$,
and thus the condition
$(ii)$ in Definition \ref{weak} becomes empty.  But this does not
satisfy the morphism property (\ref{LA_morph}) since $A_i$ clearly does
not satisfy also the second set of field equations (\ref{eqs2}) (which
would imply that all $A_i$s are closed). 
\end{example}
 Under suitable further conditions it is nevertheless possible to reverse
the above mentioned implication: In the above example the main problem
was that the given vector bundle
morphism excludes the existence of \emph{any} projectable section. 
\begin{proposition}\label{project}\cite{HigMc} If the sections of $E_1$ which are
projectable with respect to a given vector bundle morphism
$\phi :E_1\to E_2$ generate all of 
$\, \Gamma (E_1)$,  then $\phi-$relation implies the morphism property.
The assumption holds true in particular, if $\phi$ is fiberwise
surjective.
\end{proposition}
{\bf Proof.}
According to the assumption
 any
$s_1 ,s_1'\in \Gamma (E_1)$ decompose as
$s_1 =\sum a_i \xi_i , \hspace{4mm} s_1' =\sum a_i' \xi_i'$
such that $\xi_i ,\xi_i'$ are $\phi -$related to some $\eta_i ,\eta_i'\in
\Gamma (E_2)$, respectively. Since obviously (\ref{decomp}) holds
true, we should prove (\ref{Higg}):  
By $\phi-$relation  $\phi \circ [\xi_i
,\xi_j']=[\eta_i ,\eta_j']\circ \phi_0$. Using this relation in
the application of $\phi$ to 
\beq
[s_1 ,s_1']=\sum a_i a_j' ([\xi_i ,\xi_j']+
\sum \rho_1 (s_1) ( a_j') \xi_j'  - \sum
\rho_1 (s_1' )( a_j) \xi_j 
\eeq
we indeed find (\ref{Higg}).

Finally, if $\phi$ is fiberwise surjective, there exists an isomorphism
between $E_1$ and $\ker \phi \oplus \phi_0^* E_2$. Evidently, any
section of  $\ker \phi$ is $\phi-$related to the zero section of $E_2$
and all sections of $\phi_0^* E_2$ are generated by
$\phi^*(\Gamma(E_2))$; thus all sections of $\ker \phi \oplus \phi_0^*
E_2\cong E_1$ are generated by projectable sections.
$\blacksquare$


Let us notice that since the morphism equation (\ref{LA_morph}) 
and the proof above are local, the statement of Proposition
\ref{project} remains unchanged if we replace $M_1$ with an open neighborhood
of any point $x_0\in M_1$. This argument is used in the next
\begin{proposition}\label{M3} Any $\phi-$relation  with a
base map that is a local immersion is a morphism of Lie algebroids.
\end{proposition}
{\bf Proof.}
If $\phi_0\colon M_1\to M_2$ is a local immersion then for any point $x_0\in
M_1$ there exists a coordinate chart $(U, X^i)$ around $\phi_0
(x_0)$
and an open neighborhood $V\subset M_1$ of $x_0$ such that
$\phi_0 (V)\subset U$ is given by the set of equations $ X^{1+\, \dim M_1}
=\ldots =X^{\, \dim M_2}=0 $ and $\phi_0\colon V\to\phi_0 (V)$ is
a diffeomorphism.

Now one can show that any section of $(E_1)_{|V}$ is projectable with
respect to the restriction of $\phi$ on $(E_1)_{|V}$, i.e. $\phi -$related to
some section of $(E_2)_{|U}$ as a consequence of the following simple facts:
\begin{itemize}
\item The restriction of $\phi$ defines a map of sections $\Gamma
((E_1)_{|V})\to\Gamma ((E_2)_{|\phi_0 (V)})$
\item Any section of $(E_2)_{|\phi_0 (V)}$ can be extended as a section
of $(E_2)_{|U}$.
\end{itemize} $\blacksquare$

The last statement is of particular interest due to
\begin{proposition} $\phi \colon E_1 \to E_2$ is a morphism of Lie
algebroids, iff its graph $\phigra \colon E_1 \to E \equiv E_1 \boxplus
E_2$ is a morphism of Lie algebroids. 
\label{prop:gra}
\end{proposition}
{\bf Proof.}
With $\Omega_E^\cdot (M)\cong
\Omega_{E_1}^\cdot (M_1) \otimes \Omega_{E_2}^\cdot (M_2)$, $^E\md
={^E\md_1} + {^E\md_2}$, and 
$\Phigra (\omega_1\otimes\omega_2 )=\omega_1\wedge \Phi (\omega_2 )$
for all $\omega_i \in  \Omega_{E_i}^{q_i} (M_i)$ (so that
$\varepsilon_1(\omega_1) = \deg{\omega_1} = q_1$)  one has
\beqn
({{}\md_1}\Phigra -\Phigra \, {^E{}\!\md} )(\omega_1\otimes\omega_2 )=
(-1)^{q_1}
\omega_1\wedge (\md_1\Phi -\Phi\md_2)\omega_2
\eeq
which vanishes identically if and only if $\Phi$ is a chain map.
$\blacksquare$

Since the base map of $\phigra$ is even an embedding,
 the general notion of Lie algebroid morphism 
can be reduced to the simplified notion of $\phi$-relation of Lie
algebroids, $E_1  \!  \stackrel{\phigra}{\Sim}\! E$.

Finally, the chain property (\ref{LA_morph}) may be reformulated also
nicely in terms of operators living in one and the same bundle. Recall
that ${\,\!{}^E{}\Phi}$ and $^E{}\md$ both act inside $ \Omega^\cdot_E(M)$
(cf.~Eq.~(\ref{Phimaps}) and end of the previous section); while
${\,\!{}^E{}\Phi}$ is of (total) degree zero,  $^E{}\md$ is of degree one. 
We have:
\begin{proposition} $\phi \colon E_1 \to E_2$ is a morphism of Lie
algebroids, iff the induced operator ${\,\!{}^E{}\Phi}$ commutes with  $^E{}\md$
on  $ \Omega^\cdot_E(M)$, i.e.~iff the operator 
\beq {\,\!^E{}\!F_\phi} := [ \,\!^E{}\md ,{\,\!{}^E{}\Phi} ]  \label{commu} 
 \eeq
 vanishes.
\end{proposition}
{\bf Proof.}
By definition, ${\,\!{}^E{}\Phi} =P_1 \circ \Phigra$.
Since evidently  $\,\!^E{}\md \circ P_1 =  P_1 \circ \md_1$ holds
true, we obtain
\beqn  ^E\md \; {\,\!{}^E{}\Phi} - {\,\!{}^E{}\Phi} \; {}^E\md =  P_1 \circ
({{}\md_1}\Phigra -\Phigra \, {^E{}\!\md} ) \, ,  \eeq
which concludes the proof due to 
Proposition \ref{prop:gra} and the fact that $P_1$
is an injection. 
$\blacksquare$

\vskip5mm

Maybe some warning is in place: 
The above notion of a morphism, in any of its formulations,
applied to the cotangent bundle of
two Poisson manifolds, does  \emph{not} coincide with a
Poisson morphism. In contrast, a Poisson map, i.e.~a map $\hat \phi_0 \colon
(M_1,\CP_1) \to (M_2,\CP_2)$ with $(\hat \phi_0)_\ast \CP_1|_x =
\CP_2|_{\phi_0(x)} \, \forall x \in M_1$, gives rise only to a bundle
morphism $\hat \phi \colon TM_1 \to TM_2$ by means of the tangent map
$\hat \phi := (\hat \phi_0)_\ast$. This generalizes in the following
way
\begin{deff} Let $(E_i , [\cdot ,\cdot ],\rho_i )$ be Lie algebroids
over base manifolds $M_i$, $i=1,2$. We say that a bundle map
$\hat\phi\colon E_1^* \to E_2^*$ is \emph{a comorphism} if the induced
operator $\hat\Phi\colon \Gamma (E_2)\to \Gamma (E_1)$ satisfies the
following properties:
\bqa
(\hat\phi_0)^* (\rho_2 (s_2)(f))&=&\rho_1
(\hat\Phi (s_2))((\hat\phi_0)^*(f))
\qquad \forall f\in C^\infty (M_2), s_2 \in\Gamma (E_2)  \nonumber \\
\hat\Phi ([s_2,s_2'])&=&[\hat\Phi (s_2),\hat\Phi (s_2')] \qquad \qquad \;\forall
s_2,s_2'\in\Gamma (E_2) \nn
\ea
\end{deff}
In this terminology a Poisson map thus corresponds to a 
comorphism of the respective Poisson Lie algebroids, $\hat \Phi$ then
being nothing but the pullback of differential 1-forms. 
%

An algebraic generalization of these  notions may be found in \cite{Mc},
such that a morphism (comorphism) of Lie algebroids corresponds to a
comorphism (morphism) of the related pseudoalgebras, respectively.

We conclude this section with a short remark about
covariance of the
field equations (\ref{eqs1}), (\ref{eqs2}). Obviously the total set of
field equations must be covariant---they are the Euler Lagrange
equations of a completely covariant action functional, cf., e.g.,
(\ref{actiondelta}) or (\ref{actionprime}), or, likewise, they can be
reformulated frame independently as in (\ref{LA_morph}). On the other
hand, the field equations  (\ref{eqs2}) are not only not written in an
explicitly covariant form, by themselves they even are not frame
independent. The reason for this is the (kind of)
Leibniz rule satisfied by the operator  (\ref{Fphi}),
\beq \label{leibniz}
 F_{\! \phi} (\o_2 \wedge \o_2') &=& F_{\! \phi} (\o_2)  \wedge\Phi( \o_2') +
(-1)^{p} \Phi (\o_2) \wedge  F_{\! \phi} (\o_2')  \, , 
\eeq
which holds for arbitrary  $\o_2 \in \Gamma(\Lambda^p E_2^*)\, , \;  \o_2' \in
\Gamma(\Lambda^q E_2^*)$. 
Indeed, with the
abbreviations\footnote{For notation and conventions recall end of section
\ref{sec:prelim}.}
\beqa F^i &:=&  F_{\! \phi}(X^i) =  \md_1 X^i
-\rho^{i}_IA^I   \;  , \label{Fupper} \\  F^I &:=& F_{\! \phi}(b^I) = 
 \md_1 A^I + \frac{1}{2} C^I_{JK} A^J \wedge A^K \; , \label{Flower}
\eeqa
the first and second set of field equations are $F^i=0$ and $F^I =0$,
respectively.  Suppose now we
change frame from $b^I$ to a new one,  $\tilde{b}^I$,
by means of
$b^I =B^I_J \tilde{b}^J$.
Then, by means of (\ref{leibniz}), we find
$F^I = \Phi(B_J^I) \,  F^{\widetilde J} + \Phi(B^I_{J,i}) F^i \wedge
A^J$. 
This obviously implies that only upon
usage of  $F^i=0$, which itself clearly \emph{is} covariant (also
with respect to change of coordinates $X^i \to \widetilde X^i$),
we may conclude  $F^{\widetilde I}=0$ from $F^I=0$.

This may be cured by means of an auxiliary connection $\Gamma$ on
$E_2$, introducing
\beq
F^I_{\! (\Gamma)} := F^I +  \Gamma^I_{iJ} F^i
\wedge A^J \, . \label{FGamma}
\eeq
This option shall be investigated into further depth in a
separate paper \cite{AT1}. In the present paper we are
interested particularly in morphisms, $F^i=0=F^I$, in which case
covariance of (\ref{Flower}) is of subordinate importance. The issue
of covariance will become more important in the context of
the following section, however.


%




\section{Generalized gauge symmetries}
\label{sec:Generalized}
We now turn to interpreting and generalizing
the gauge symmetries of the
PSM. In view of the generalization (\ref{Fupper}) and
(\ref{Flower}) of the field equations (\ref{eqs1}) and (\ref{eqs2}),
it is suggestive to replace the gauge symmetries  (\ref{sym1}) and
 (\ref{sym2}) by
\beq
\de^0_\e X^i &=& \rho^i_I \e^I \, ,\label{Sym01}\\
\de^0_\e A^I &=&  \md_1\e^I + C^{I}_{JK} A^J \e^K \,\label{Sym02} ;
\eeq
without further mention it is assumed furthermore that $\de^0_\e$
obeys an (ungraded) Leibniz rule (which is used e.g.~when establishing
gauge invariance of (\ref{action}) up to boundary terms).

Such as we were able to cast 
(\ref{Fupper}) and (\ref{Flower}) into a more elegant and
covariant form, cf., e.g., (\ref{commu}), 
 and prove the equivalence of their vanishing with the
morphism property of Lie algebroids, we may now strive for similar
issues in the context of  (\ref{Sym01}) and  (\ref{Sym02}). This indeed
is part of the intention of the present section. However, first we need
to notice that in the context of symmetries the non-covariance of the
formulas (\ref{Sym01}),  (\ref{Sym02}) or  (\ref{sym1}),
 (\ref{sym2}) is much more severe than in the case of the field
 equations, which are only not written in explicitly covariant form in
 (\ref{eqs1}), (\ref{eqs2}), while, as a total set, they certainly are
 covariant. As written, the symmetries either
 have only on-shell meaning (when there is an action functional like
 in the PSM this is tantamount to having meaning
 only as quotient of all symmetries modulo so-called trivial ones,
cf.~also \cite{PSMRiemann}) or they are defined only for trivial or
flat bundles
 $E_2$ (respectively for topologically rather trivial Poisson manifolds)!

Let us be more explicit about this: An infinitesimal gauge symmetry such
as  (\ref{Sym01}) and (\ref{Sym02}) is supposed to be a vector field on
the (infinite dimensional) space $\cM = \{\phi \colon E_1 \to E_2\}
\cong \{ \Phi \}$ of fields and thus, for a fixed element
 $\phi$  in $\cM$, a vector $\cV \in T_\phi \cM$. Note that $\cM$ is a
bundle over $\cM_0 = \{\phi_0 \colon M_1 \to M_2\}$, the space of base
maps. The projection of $\cV$ to $\cM_0$ then gives a vector 
$\cV_0 \in T_{\phi_0}\cM_0$. Eq.~(\ref{Sym01})
indeed corresponds to a vector on  $\cM_0$, as may be seen by changing
coordinates on $M_2$ (or likewise also local frames in $E_2$). However, 
 (\ref{Sym01}) and (\ref{Sym02}) together
 do \emph{not} give a well-defined vector on
the total space $\CM$. Indeed, if we change frame in $E^*_2$, $b^I=B^I_J
\tilde {b}^J$, such that $\e^I = B^I_J(X(x)) \,
\tilde {\e}^J$ etc, a straightforward calculation yields
\beq \md_1\e^I + C^{I}_{JK} A^J \e^K = B^I_J \left( \md_1\tilde \e^J +
\tilde C^{J}_{KL} \tilde A^K \tilde \e^L \right) + B^I_{J,i} \md_1X^i
\tilde {\e}^J + B^I_{J,i} \tilde \rho^i_K
\tilde A^J  \tilde  \e^K - B^I_{J,i} \tilde \rho^i_K
\tilde A^K  \tilde \e^J \, ; \label{rhs}
\eeq
on the other hand, by the Leibniz rule we obtain,
\beq \d^0_\e (B^I_J \,
\tilde {A}^J ) =B^I_J \,
\d^0_\e \tilde {A}^J + B^I_{J,i} \tilde {A}^J
\d^0_\e X^i \, . \label{Leibniz0}
\eeq
The difference of the right hand sides of (\ref{rhs}) and
(\ref{Leibniz0}) is
\beq \label{difference}
B^I_{J,i} F^i
\tilde {\e}^J \, .
\eeq
Therefore, in general (\ref{Sym01}) and (\ref{Sym02}) do not provide a vector
in $T_{\phi_0}\cM_0$; it is globally well-defined only on fields
satisfying $F^i=0$ or when  $B^I_J$ can be  chosen
consistently to be  $X$-independent. The first option is (part of) the on-shell
condition, the second one corresponds to the existence of a flat
connection in $E_2$. 
In this case (\ref{Sym02}) depends implicitly on the frame
and on the flat connection chosen, which is zero in the particular
frame chosen, but becomes non-zero if we change the frame. 

At this point let us emphasize that $\delta_\epsilon$
is \emph{not} a tangent vector field to $\cM$
if it satisfies $\delta_\epsilon A^I
= B^I_J\delta_\epsilon \tilde A^J$ (which would correspond to the
absence of all three terms in (\ref{rhs}))
with respect to a change of frame 
 $b^I=B^I_J
\tilde {b}^J$; it is an element of $T_\phi \cM$ only when it satisfies
an ungraded Leibniz rule, i.e.~in particular
\beq \delta_\epsilon A^I =  B^I_J\delta_\epsilon \tilde A^J+
 B^I_{J,i} \tilde {A}^J
\d_\e X^i \,  \label{Leibniz1}
\eeq
 (which would correspond to the
absence of the last and the third to last term in (\ref{rhs}), which
together combined into (\ref{difference})). As a consequence, even if
one uses a connection on $E_2$ to provide a global and frame
independent definition of the tangent vectors $\delta_\e$, the
explicit formula for $\delta_\e A^I$ will \emph{not} be covariant
(in the usual sense) with respect to capital indices (containing only
covariant derivatives and $E_2$-tensors).\footnote{There is one 
trivial exception to this statement, namely the case for which the
second term in (\ref{Leibniz1}) vanishes identically (for all choices
of $B^I_J(X)$): This happens iff $\delta_\epsilon X^i \equiv 0$ for
all $\e$, which, in view of the covariance and off-shell validity of
(\ref{Sym01}), in turn is tantamount to $\rho \equiv 0$, i.e.~this
happens iff $E_2$ is a bundle of Lie algebras.} In contrast,  $\delta_\e
X^i$ is covariant with respect to $i$, since multiplication by (the
pullback of) the Jacobian
 of a coordinate change on $M_2$ is in 
agreement with the Leibniz property of $\delta_\e$.

For the rest of the section,
we will proceed as follows: 
In view
of the above observation, $\delta^0_\e$ as defined in (\ref{Sym01}) and
(\ref{Sym02}) should have a good, more abstract
\emph{on-shell} interpretation. Indeed, we will see that it
corresponds to an infinititesimal homotopy of Lie algebroid
morphisms. Simultaneously this picture provides an on-shell
integration of the infinitesimal symmetries $\delta^0_\e$.
Next we  want to lift the on-shell symmetry to a
well-defined off-shell symmetry. This is not unique certainly. One
option is to do this in such a way that the (infinitesimal)
inner automorphisms of $E_1$ and $E_2$ are contained as Lie subalgebras. 
This will turn out to be done most efficiently in terms of $E$-Lie
derivatives of the exterior sum Lie algebroid $E=E_1 \boxplus E_2$. The second
option is to employ a connection on $E_2$, such that 
for flat connections $\Gamma$, and in a frame for which $\Gamma =0$,
one reobtains the original formulas for $\delta^0_\e$. This second
option shall be mentioned at the end of this section 
peripherically only; for more details we refer to  \cite{AT1}.

\vskip7mm


\begin{deff}
Let $E_1$ and $E_2$ be Lie algebroids  over smooth manifolds
$M_1$ and $M_2$,
respectively. We say that the two
morphisms $\phi ,\phi' \colon E_1\to E_2$ 
are \emph{homotopic},
 iff there is  a morphism $\overline\phi$ from the  Lie algebroid
  $\bE := E_1 \boxplus TI$ over the manifold $N=M_1\times I$,
  $I \equiv[0,1]$, such that
the restriction of $\bphi $ to the boundary components $M_1\times\{0\}$ and
$M_1\times\{1\}$ coincides with
$\phi$ and $\phi'$, respectively. 
\end{deff}


\begin{proposition} \label{prop:global}
Two Lie algebroid morphisms 
 $\phi $ and $\phi '$
are homotopic, iff they can be connected by a flow of $\delta^0_\cdot$ as
defined in (\ref{Sym01}),  (\ref{Sym02}). 
\end{proposition}
Note that, as outlined above, $\delta^0_\e$ is well-defined on-shell,
i.e.~as a vector field
on the subset of $\CM$ satisfying the field equations $F^i=0=F^I$; in
the above proposition $\delta^0_\e$ is understood in this on-shell sense. 
\vskip2mm \noindent
{\bf Proof.}
Given a local frame $\{b_I\}$ in $E_2$ over a coordinate chart $\{X^i\}$,
we immediately obtain the following system of equalities from the
chain property of $\bar \phi$:
\beq \label{bareom}
\bF^i
= \,\!^{\overline E}\md X^i - \rho_{I}^i \bA^I \equiv 0, \hspace{4mm}
\bF^I = \,\!^{\overline E}\md\bA^I + \frac{1}{2}C_{JK}^I\bar
A^J\wedge\bA^K
\equiv 0 \, ,
\eeq
where the structure functions $C_{JK}^I$ and $\rho_{I}^i$ depend on
$X(x)$, $x \in M_1$, but not on $t$.  On the other hand, 
by definition, 
$\bE = E_1\boxplus TI$ and 
${\,\!^{\overline E}\md}_U =\md_1 +\md t\wedge\pt_t$; correspondingly, 
$\overline A^I = A^I + \overline A_t^I \md t$, with  $A^I \equiv
A^I(t) \equiv A_\alpha^I b^\a$ being (local) $t$-dependent
$E_1$-1-forms. {}Adapting (\ref{bareom}) to this splitting, and
renaming  $\overline A_t^I$ to $\e^I$, 
 we obtain
\beq\label{bar_eq1}
 \bF^i &=& F^i(t) + 
\md t \left( \pt_t X^i -\rho_{I}^i \e^I_t\right) \\\label{bar_eq2}
\bF^I &=& F^I(t) +
\md t\wedge \left( \pt_t A^I-\md_1 \e^I +
C_{JK}^I\e^J\bA^K\right) 
\eeq
where $F^i$ and $F^I$ are of the form (\ref{Fupper}) and (\ref{Flower})
and $ \pt_t A^I \equiv (\pt_t \bA^I_\a ) b^\a$.
This proves that $\bF^i = 0=\bF^I$, iff for any $t$ one has
$F^i = 0=F^I$ and
$\pt_t X^k = \de^0_\e X^k$, $\pt_t A^K =
\de^0_\e A^K$.
$\blacksquare$

If $M_i$ are manifolds with boundary one has to take
care about boundary conditions. In particular, the space of morphisms
from $TI$ to an arbitrary Lie algebroid $E$ over a manifold $M$ modulo
homotopies (with fixed boundary contribution) gives the {\sl fundamental}
or {\sl Weinstein's groupoid} of $E$, cf.~\cite{CFHam}. 
Thus, 
the on-shell part of gauge symmetries (\ref{Sym01}), (\ref{Sym02})
is well-motivated now. It corresponds to the infinitesimal flow of a
homotopy of Lie algebroid morphisms.

\vskip7mm

We now turn to a possible off-shell definition of the
gauge symmetries without the introduction of any further structures
such as a connection on $E_2$, employed in an alternative approach in
\cite{AT1}. 
Concretely this means that we
want to  extend (\ref{Sym01}), (\ref{Sym02}) to a differential
$\delta_\e$, satisfying (\ref{Leibniz1}), where for $F^i=0=F^I$ 
the gauge transformation
$\delta_\e$ reduces to $\delta^0_\e$---and we want to relate this
differential on field space to a differential operator on or between
finite dimensional bundles, in analogy of what we did with the field
equations. 

\begin{deff} 
We call an operator
$\cV \colon \O^\cdot_{E_2}(M_2) \to  \O^{\cdot +
\deg{\cV}}_{E_1}(M_1)$ a \emph{$\Phi$-Leibniz} operator, if it satisfies
$\forall \o , \o' \in \O_{E_1}^\cdot(M_1)$ ($\o$ homogeneous) 
\beq \cV (\o \wedge \o') = \cV (\o) \wedge \Phi(\o') + (-1)^{\deg{\cV}
\deg{\o}} \;
\Phi( \o) \wedge \cV (\o')  \, , 
 \eeq
and likewise an operator $\EcV$ in
$\Omega^\cdot_E(M)$ (of 
fixed degree)  \emph{$\EPhi$-Leibniz}
if it satisfies the above equation with $\cV$ and $\Phi$ replaced by
 $\EcV$ and $\EPhi$, respectively. 
\end{deff}
An example for a degree one $\Phi$-Leibniz operator is
provided by $F_\phi$, cf.~Eq.~(\ref{leibniz}); likewise
${\,\!^E{}\!F_\phi}$, defined in ({\ref{commu}), is
$\EPhi$-Leibniz. More generally, obviously any consecutive 
 application  (in both possible orders)
 of a (standard) Leibniz operator with $\Phi$ ($\EPhi$) 
 gives a $\Phi$-Leibniz ($\EPhi$-Leibniz) operator.
\begin{deff} \label{def:gauge}
We call $\delta \Phi \colon
\O^\cdot_{E_2}(M_2) \to  \O^{\cdot}_{E_1}(M_1)$ an \emph{infinitesimal
gauge symmetry}, if it is a degree zero $\Phi$-Leibniz operator satisfying
$\md_1 \delta \Phi \approx \delta \Phi \md_2$, where $\approx$ denotes
an on-shell equality (i.e.~it has to be an equality for all $\Phi$
with $F_\phi = 0$). Likewise a degree zero $\EPhi$-Leibniz operator
$\EdPhi$ is
an infinitesimal gauge symmetry if it satisfies 
\beq
[\EdPhi, \Emd ] \approx 0 & \Leftrightarrow & [\EdPhi, \Emd ]|_{\phi
\colon\!\! [\EPhi,\Emd]=0} =0 \label{eq:gauge}
\eeq
and  $\mathrm{im}\EdPhi \subset \mathrm{im} P_1$,
$\EdPhi \circ P_1 =0$.
\end{deff}
This is motivated as follows: $\delta \Phi
\sim \md \Phi_t/\md t|_{t=0}$ for some family of $\Phi$s parametrized
by $t$. Correspondingly, since $\Phi$ is of degree zero, also $\delta
\Phi$ is, and functoriality of $\Phi$, $\Phi(\o \wedge
\o') =\Phi(\o) \wedge \Phi(\o')$, results in the
$\Phi$-Leibniz property. Finally,
$\Phi_t$ satisfying the field equations implies that  $\delta \Phi$
does so on use of the field equation for $\Phi \sim \Phi_{t=0}$.
All this applies analogously 
to $\EdPhi$, where, however, in
addition we need to take care of the fact that $\EPhi$ is not an
arbitrary operator in $\OE$, but restricted as specified in
(\ref{PhiE}) and the text thereafter. 

One of the main features of a gauge symmetry is that it maps solutions
of field equations into solutions. Here, the solutions have the
meaning of a morphism (of Lie algebroids) $\phi \colon E_1 \to E_2$.
To construct gauge symmetries we may thus proceed as follows: Let
the gauge transformed morphism $\tilde \phi$ be given by  $\tilde \phi
:= (a_1)^{-1} \circ \phi \circ a_2$ where $a_i \in \mathrm{Aut}(E_i)$,
$i=1,2$, the respective group of automorphisms of $E_i$. This defines
a right action of $\mathrm{Aut}(E_1) \times \mathrm{Aut}(E_2)$ on $\cM
= \{ \phi \}$, which on the level of Lie algebras provides a
\emph{homomorphism} $\mathrm{aut}(E_1) \oplus \mathrm{aut}(E_2) \to
\Gamma(T\cM)$.

A subgroup of the automorphism group of a Lie algebra
$E_i \cong \g_i$ is the group of  inner automorphisms, given by the adjoint
action of the Lie group $G_i$ which integrates $\g_i$;
infinitesimally, this is just the regular representation of the Lie
algebra $\g_i$, i.e.~the action of $\g_i$ onto itself
given by multiplication in the Lie algebra, $v_i \mapsto [  v_i , \cdot]$ (a homomorphism
of $\g_i \to \mathrm{aut}(\g_i)$).
Although not every Lie algebroid has a (sufficiently smooth) Lie
groupoid integrating it
(cf.~\cite{CFm1} for the necessary and sufficient
conditions), we still may generalize the infinitesimal picture to the
setting of Lie algebroids:
Given a section $s_i \in \Gamma(E_i)$, we may regard $\ELi_{s_i}$ as a
vector field on $E_i$, which, due to $\ELi_{s_i}(s_i') = [s_i,s_i']$
and the Jacobi property of the Lie algebroid bracket, is an
infinitesimal automorphism of $E_i$.

That $\ELi_{s_i}$ indeed can be
regarded as a vector field on $E_i$ may be seen as follows:
$C^\infty(M_i)$ and $\O^1_{E_i}(M_i)$ are  fiberwise constant and
linear functions on $E_i$, respectively. Together they
generate all of $C^\infty(E_i)$. Local coordinates $X$ on $M_i$
and a local coframe $b^I$ provide a local coordinate system on $E_i$.
Applying a vector field to local coordinates gives its components in
this coordinate system; these components may be easily extracted from 
Eq.~(\ref{Lie_2}), showing that they are linear in the fiber
coordinates. The $E_i$-Lie derivative $\ELi_{s_i}$
provides a uniquely defined
lift of $\rho(s_i) \in \Gamma(TM_i)$ to  $\Gamma(T(E_i))$; in contrast
to the lift given by a contravariant connection this lift is not
$C^\infty$-linear in $s_i$, certainly.

\begin{proposition} \label{SymmProp}
For arbitrary sections $s_i \in \Gamma(E_i)$, $i=1,2$, 
\beq \label{sonetwo}
 \delta \Phi :=  \Phi \circ \ELt_{s_2} - \ELo_{s_1} \circ \Phi 
\eeq
is an infinitesimal gauge symmetry. For any $\Phi \in \cM$, its action
on a local coordinate system  $X^i,b^I$ on $E_2$ defines a Leibniz
operator $\d_\e$ (an element in $\Gamma(T_\phi \cM)$), which 
agrees with $\d^0_\e$  given in
(\ref{Sym01}), (\ref{Sym02}) on-shell where $\e = \Phi^!(s_2) - \iota_{s_1} A$
(and $A \equiv \Phi^!(\delta)=A^I \otimes \barb_I$). Moreover, the
commutator of two such infinitesimal gauge transformations is again of
the same form, $[\d_\e , \d_{\e'} ] = \d_{\e''}$, where $\e''$ results
{}from $s_1'' = [s_1,s_1']$ and  $s_2'' = [s_2,s_2']$. 
\end{proposition}
The statement in this proposition may be simplified by saying that
there exists a \emph{homomorphism} 
$\,\Gamma(E_1) \oplus \Gamma(E_2) \to \Gamma(T_\phi \cM),\, \delta_{(s_1,s_2)}
\Phi \mapsto \delta_\e$; however, we refrained from doing so, since,
at least at this point, we
did not want to go into the details of defining properly the infinite
dimensional tangent vector bundle $T \cM$ (while still we will come
back to this perspective in more detail below). 
Let us remark already at this point, moreover, that the set of $\e$'s
that one may obtain in this fashion is too restrictive, yet. Assume
e.g.~that $\phi$ corresponds to $A^I=0$ and $X^i(x)=const$. Then any
$\e$ of the above form is necessarily constant, while it need not be
so in (\ref{Sym01}), (\ref{Sym02}),  where $\e \in
\Gamma(M_1,\phi_0^* E_2)$ arbitrary.

\proof
First it is easy to see that (\ref{sonetwo}) provides an infinitesimal
gauge symmetry in agreement with definition \ref{def:gauge}. As a
composition of Leibniz operators with $\Phi$ it is $\Phi$-Leibniz, and
since $\ELi$-Lie derivatives commute with the respective differential
$\md_i$, $\md_1 \, \Phi \approx \Phi \,\md_2$ is seen to result in 
$\md_1  \,\delta \Phi \approx \delta \Phi \, \md_2$.

To determine the desired map $(s_1,s_2) \in \Gamma(E_1) \oplus \Gamma(E_2)$ to
$\e \in \Gamma(M_1,\phi_0^* E_2)$, we may use Cartan's magic formula 
(\ref{Lie_derivative}) to rewrite $\delta \Phi= \delta_{(s_1,s_2)}
\Phi$ according to
\beq \label{delta}
\delta_{(s_1,s_2)}
\Phi = \Phi  \, \ELt_{s_2} - \ELo_{s_1} \, \Phi =  \delta^0_{(s_1,s_2)}
\Phi - \left(
F_{\! \phi} \,
\iota_{s_2} + \iota_{s_1} \, F_{\! \phi}\right) \, ,
\eeq
where
\beq \label{delta0}
\delta^0_{(s_1,s_2)}\Phi \equiv \md_1 \left(\Phi
\, \iota_{s_2} - \iota_{s_1} \, \Phi \right) +\left(\Phi \, 
\iota_{s_2} - \iota_{s_1} \, \Phi \right)  \md_2 \, .
\eeq 
While the last two terms in (\ref{delta}) vanish on-shell obviously,
it is easy to verify that $\delta^0_{(s_1,s_2)}$ acting on $X^i$ and
$b^I$ agrees with  $\delta^0_\e$ in 
(\ref{Sym01}), (\ref{Sym02}) with the parameter $\e$ as given above.
Finally, since actions coming from the right and from the left
commute, it is obvious that $[\delta_\e , \delta_{\e'}]$ (with 
 $\e$ and $\e'$ of the given form) when applied
to $A^I$ and $\phi_0^* X^i$ is tantamount to the application of  
$\Phi \circ [\ELt_{s_2},\ELt_{s_2'}]  - [\ELo_{s_1},\ELo_{s_1'}]
\circ \Phi$ to $b^I$ and $X^i$, respectively. The statement now
follows since $E_i$-Lie derivatives are a representation of
$\Gamma(E_i)$. 
$\blacksquare$

Note that in contrast to $\delta^0_\e$, the operator
$\delta^0_{(s_1,s_2)}\Phi$ in Eq.~(\ref{delta0})
is defined frame-independently. However, it
now is not a $\Phi$-Leibniz operator (only on-shell it is).
We remark in parenthesis that
one may also generalize the operator in Eq.~(\ref{delta0}) to one
defined for arbitrary sections $\e \in \Gamma(M_1,\phi_0^* E_2)$:
$\delta_\e^0\Phi := \md_1 \, i_{\e} + i_\e \, \md_2$ with the
operator $i_\e$ being defined by means of
\beq \label{iota}
i_\e(f b^{I_1}  \wedge \ldots \wedge  b^{I_k}) :=
\sum\limits_{j=1}^k  (-1)^{j+1} \e^{I_j} \, \Phi \left( f b^{I_1} \wedge
\ldots
\wedge \widehat b^{I_j} \wedge \ldots \wedge  b^{I_k}  \right) \, .
\eeq
But such somewhat artificial constructions do not seem very
promising. Instead, the right step is to take recourse 
to the exterior sum bundle $E = E_1 \boxplus E_2$. This has the effect that in
the end the section $\Phi^!s_2|_x=s_2(X(x))$, $x\in M_1$, of the
previous proposition is replaced
by a likewise section that depends  on both variables, $X(x)$ and
$x$, independently.

\begin{theorem} \label{theorem1}
Any section $\e \in \Gamma(E)$ which is projectable to a section
of $E_1$ ($p_1$-projectable) defines an infinitesimal gauge symmetry
by means of 
\beq
\EdePhi := [\EPhi,\EL_\e ] \, , 
\eeq
and the commutator of two such gauge transformations for $\e$, $\e'$
is the gauge
transformation associated to $[\e,\e'] \in \Gamma(E)$. In particular, for
``vertical'' sections $\e \in
\Gamma(\mathrm{pr}_2^*E_2) \subset \Gamma(E)$ its action on local
fields $X^i$, $A^I$ is given by:
\beq \label{newsymm} \delta_\e X^i = \delta^0_\e X^i
\quad , \qquad
 \delta_\e A^I = \delta^0_\e A^I - \e,_i F^i \,, \eeq where $\delta^0$ was
 defined in Eqs.~(\ref{Sym01}),  (\ref{Sym02}) and $F^i \approx 0$ in
 Eq.~(\ref{Fupper}). 
\end{theorem}

\proof Obviously $\EdePhi$  is $\EPhi$-Leibniz, and it obeys
Eq.~(\ref{eq:gauge}) since  $\Emd$ commutes with any $E$-Lie derivative
and on-shell (by definition) also with $\EPhi$. Thus it remains to
check the final two restrictions on an infinitesimal gauge
transformation  specified in definition \ref{def:gauge}. It is these
conditions that make the restriction to $p_1$-projectability (as
defined in the beginning of Sec.~\ref{sec:morph}, where the bundle
map $\phi$ is
replaced by $p_1 \colon E \to E_1$, cf.~diagram 1) of $\e \in
\Gamma(E)$ necessary. To see this we first split $\e$ according to
$E = \mathrm{pr}_1^*E_1 \oplus  \mathrm{pr}_2^*E_2$ into $\e=\e_1 +
\e_2$ and use linearity in $\e$. Due to $ [\EPhi,\EL_{e_2}]= \EPhi\EL_{e_2}$,
the image of $\EdetPhi$ lies trivially in $\mathrm{im}P_1= \mathrm{im}\EPhi$,
and also obviously it acts trivially on $P_1(\o_1) = \o_1 \otimes 1$
for all $\o_1 \in \OEo$. To ensure that also  $\EdeoPhi$ kills all
$\o_1 \otimes 1$, we introduced the commutator of $\EL_{\e_1}$ with $\EPhi$,
the latter operator acting as the identity on the image of
$P_1$. However, in this case both conditions are satisfied if and only
if  $\e_1$ depends on $x \in M_1$ only, but not also on $X \in
M_2$ (consider e.g.~$\EL_{\e_1} \EPhi = \i_{\e_1} \md   \EPhi +
\ldots$); more abstractly this means that $\e$ is $p_1$-projectable,
the corresponding $E$-Lie derivative generating only automorphisms of
$E$ that are preserving fibers over $M_1$.

Two successive gauge transformations with parameter $\e$ and $\e'$ are
characterized by the operator 
$[[\EPhi,\EL_\e ],\EL_{\e'} ]$.\footnote{That the successive
application of a vector field in field space $\cM$ has again such a simple
operator-description (being a second order differential operator on $\CM$, it
now is no more $\EPhi$-Leibniz, certainly, but satisfies a similar
higher analog of this property), is also a benefit of the present
approach using operators on $\OE$, $E=E_1 \boxplus E_2$.}
Subtracting from this the
corresponding operator with  $\e$ and $\e'$ exchanged and using the
Jacobi condition for the (graded) commutator bracket, we obtain 
$[\EPhi,[\EL_\e ,\EL_{\e'}]]= [\EPhi,\EL_{[\e,\e']} ]$,
a gauge transformation with parameter $[\e,\e']$.

To relate the gauge transformations above to explicit transformations
acting on the fields, we proceed similarly to before
(cf.~Eqs.~(\ref{delta}) and (\ref{delta0})), where now the
splitting becomes a bit more elegant:
\beq
\EdePhi &=& [\EPhi,[\Emd,\i_\e ]] = \EdzePhi-[\EF,\i_\e]  \, ,  \\
 \EdzePhi &\equiv&  [\Emd,[\EPhi,\i_\e] ] \, , 
\eeq
where we made use of the (graded) Jacobi property and the definition
(\ref{commu}) for $\EF$. Upon action on $X^i$, $b^I$ (or, more
generally, the image of
$P_2\colon \OEt \to \OE$)---and for $\e = s_1 + s_2$---the operator 
$\EdzePhi$ is identified easily with the one in (\ref{delta0}); for
general $p_1$-projectable $\e$ it just provides formulas 
 (\ref{Sym01}) and  (\ref{Sym02}). The on-shell vanishing
 contributions, necessary to render the gauge transformation globally
 defined  and Leibniz, are now easily calculated to be
 \beq
 [\EF,\i_\e] X^i &=&  \i_{\e_1} F^i \, ,\label{fullgauge1}\\
 {}[\EF,\i_\e] b^I &=& -\e_2^I,_i F^i 
 + \i_{\e_1} F^I \, .\label{fullgauge2}
 \eeq
Note that here we used that $\EF (\e_2^I)$ contributes only by its
derivative with respect to $X$, but not also with respect to $x$; the
latter terms cancel in the commutator
(\ref{commu}). Eq.~(\ref{newsymm}) now follows by specialization to
$\e = \e_2$. 
$\blacksquare$

Given sections $s_i \in \Gamma(E_i)$, $i=1,2$, there is a natural
inclusion as sections of the exterior sum $E$ of $E_1$ and $E_2$. With
$\e: = s_1 + s_2$ it is easy to see that the action of $\EdePhi$ on $X^i$
and $b^I$ precisely reduces to $\delta \Phi$ as given in
(\ref{sonetwo}). The extension of the present approach is that
now $s_2$ may effectively depend also on $x$ (and that due to using
the graph both, the action from the left and the action from the right
in Prop.\ \ref{SymmProp} now come from the
right); due to this $x$-dependence of $\e_2$ (while $\e_1$
is still not permitted to depend on $X$), the total action is no
more a direct sum of $\Gamma(E_1)$ with $\Gamma(E_2)$ as in  Prop.\
\ref{SymmProp}, but a semidirect sum, spanned by the two Lie
subalgebras generated by $\e_1$ and $\e_2$, respectively. 

It is needless to say that an explicit verification of the closure 
of the symmetries (\ref{newsymm}) (or even as the one with
$\e=\e_1+\e_2$, cf.~Eqs.~(\ref{fullgauge1}), (\ref{fullgauge2})) would
be a formidable task. This now was reduced to a simple line only. We
may even use the above approach to simplify the likewise calculation
of the commutator of the
initial symmetries (say in a flat bundle or used in one
particular coordinate patch):
\vskip2mm 
\noindent {\bf Corollary:} 
The commutator of  two symmetries (\ref{Sym01}), (\ref{Sym02}) 
corresponding to $\e, \e' \in \Gamma(\phi_0^* E_2)$ is 
\beq
[\d^0_\e , \d^0_{\e'} ] X^i = \d^0_{[\e,\e']} X^i \; , \quad [\d^0_\e ,
\d^0_{\e'} ] A^I = \d^0_{[\e,\e']} A^I - C^I_{JK,i} F^i \e^J {\e'}^K ,
\eeq
where $[\e,\e']^I:= \phi_0^* \left( C^I_{JK}\right)\e^J {\e'}^K$.
\proof Any section  of $\Gamma(\phi_0^* E_2)$ can be regarded as the
restriction of some section in $\Gamma (M, \mathrm{pr}_2^* E_2)$ to
the graph of $\phi_0 \colon M_1 \to M_2$ inside $M$. Notice that this
choice is not unique, certainly; given a flat connection on $E_2$, or
in a particular local frame $b^I$, (which underlies the definition of
$\delta^0$!), we can choose this extension to be constant along
$M_1$-fibers or independent of $X$. We denote these extensions again
by the same letters. Note that the bracket in $E$ induces the bracket
as specified above when restricted to the $\phigra_0(M_1)\subset M$;
however, the  bracket $[\e,\e'] \subset \Gamma(M, \mathrm{pr}_2^*
E_2)$ is in general \emph{not} constant along $M_1$-fibers; in general
it depends on $X$ due to the $X$-dependence of the structure functions
$C^I_{JK}$. By use of Eq.~(\ref{newsymm}) we thus obtain immediately
\begin{eqnarray}
 [ \de^0_\e , \de^0_{\e'} ] X^i = [\de_\e , \de_{\e'}]X^i = \de_{[\e
 ,\e']}X^i =\de^0_{[\e ,\e']}X^i \; ,
 \\  {}[\de^0_\e , \de^0_{\e' }]A^I
 =
 [\de_\e , \de_{\e'}]A^I = \de_{[\e ,\e']}A^I =
 \de^0_{[\e ,\e']}A^I - C^I_{JK,i} F^i \e^J \e^K \; .
\end{eqnarray}
$\blacksquare$

In the particular case of the PSM this reproduces the well-known
contribution rendering the algebra to be an ``open'' algebra. We now
see that this may be avoided by the additional contribution in
(\ref{newsymm}) at the cost of keeping track of the $X(x)$ dependence
of $\e$, which, however, anyway cannot be avoided in the case of a
general, non-flat bundle $E_2$.


\vskip7mm

Summing up, we  see that the gauge symmetries (\ref{newsymm}) are
well-defined off-shell and globally. They are one possible off-shell
extension of the always defined on-shell version, recognized above as a
homotopy. Another extension is provided by a connection on
$E_2$. In rather explicit terms this takes the form (besides the
obvious $\delta^{(\Gamma)}_\e X^i  =  \delta_\e X^i$):  
\beq \delta^{(\Gamma)}_\e A^I = \delta^0_\e A^I + \Gamma^I_{iJ} F^i
\e^J \, . \label{deltaGamma} \eeq
Let us remark that similarly to our considerations about
homotopy---but without requiring $F$ to vanish---it is possible to view
these transformations as the components of the covariant curvatures
$F^i$ and $F^I_{\! (\Gamma)}$ in a $(1+\dim (M_1))$--dimensional
spacetime, cf.~Eq.~(\ref{FGamma}). For a more detailed and coordinate
independent explanation of this alternative we refer to \cite{AT1}.

For both off-shell extensions it is clear by construction that they
map solutions to the field equations into other solutions. However, it
is not clear that, when specialized to the PSM, they would leave
invariant the action functional (since then the invariance needs to hold
off-shell). In fact, if e.g.~one wants to check invariance of the PSM
action (\ref{action}) with respect to (\ref{deltaGamma}), specialized
to the Poisson case, one finds invariance for all $\e_i = \e_i(x,X(x))$ if
and only if the connection $\Gamma$ is torsion-free. 

We now want to discuss the same issue for the case of
(\ref{newsymm}), also in a more coordinate independent way. For this
purpose we return to  (\ref{actiondelta}), rewriting it, however, in a
way more suitable to the graph map $\phigra$ (we prefer to use
$\phigra$  here instead of $\Ephi$, since for an action functional
we need a volume form on $M_1$, not a form on all of $M=M_1 \times
M_2$). We first remark that the joint map $\Alt \circ \Phi^*$ can be
obtained also
as the dual map to $\widetilde \phi := \phi \oplus (\phi_0)_\ast \colon
TM_1 \to T^* M_2 \oplus TM_2$. Indeed, the induced map $\widetilde \Phi$ 
then just maps $\Gamma(\Lambda^\cdot(TM_2)) \otimes \Omega^\cdot(M_2)$
to forms over $M_1$ and  $\widetilde \Phi = \Alt \, \Phi^*$. Next, we
may repeat the steps above for the map $\phigra$ instead of $\phi$
by replacing the
target Lie algebroid in the map $\phi \colon E_1 \to E_2$
by $E = E_1 \boxplus E_2$. So, $\widetilde {\phigra} = \phigra
\oplus (\phigra_0)_\ast \colon TM_1 \to E \oplus TM$
and  $\widetilde \Phigra$ 
acts from $\Gamma (\Lambda^\cdot (E\oplus
TM)^*)\cong  \Omega^\cdot_E (M)\otimes \Omega^\cdot (M)$
 to $\Omega^\cdot (M_1)$. In this way we obtain
\beq
 S[\phi] = \int_\S \widetilde \Phi (\delta+\cP)
 = \int_\S \widetilde \Phigra (\delta+\cP)
 \, . \label{actiondelta2}
\eeq
To determine the variation of  $\widetilde \Phigra$ with respect to
a gauge transformation, we first need to extend the $E$-Lie derivative
$\EL$ defined on $E$ to $E \oplus TM$ (which is \emph{not} a
Lie algebroid itself in general), i.e.~to define $\widetilde \EL_\e$
 on elements of $\Omega^\cdot_E (M)\otimes \Omega^\cdot (M)$ for any
$\e \in \Gamma(E)$: let $\widetilde \EL_\e$ restrict to $\EL_\e$ on 
 $\Omega^\cdot_E (M)$ and act as $L_{\rho(\e)}$ on $\Omega^\cdot (M)$;
 this gives a well-defined action on the tensor product since the two
 actions agree on functions. Then for any projectable 
section $\e\in \Gamma (E)$ one has 
$\de_\e(\widetilde \Phigra) =\widetilde \Phigra\, \EL_\e -
L_{\mathrm(p_1)_* (\e)} \, \widetilde \Phigra$,
where $\mathrm(p_1)_* (\e)\in
\Gamma (TM_1)$  is the projection of $\e$ to $E_1=TM_1$, and $L$
denotes the ordinary Lie derivative. The second contribution in 
$\de_\e(\widetilde \Phigra)$ takes care of the fact that one
respects the graph property.  Now we are ready to state
\begin{proposition} The PSM action  (\ref{actiondelta}) or
(\ref{actiondelta2}) is invariant
with respect to the gauge transformations (\ref{newsymm}), if the
projectable section $\e \in \Gamma(E) \cong \Omega^1(M_2)$
satisfies \beq\label{gauge_condition}
\left(\CP^\sharp \otimes \mathrm{id}  \right) \; (\md_2\e) =0\, ,
\eeq
where $\md_2$ is the de Rham operator over $M_2$ (extended trivially
to $M = M_1 \times M_2$).
\end{proposition}
\proof
In this situation we now have the identifications: $E_1 =TM_1$,
$E_2=T^*M_2$, $E=TM_1\boxplus T^*M_2$, and $\Omega^m_E (M)
=\oplus_{p+q=m} \Omega^p (M_1)\boxtimes \Gamma (\Lambda^q TM_2)$.
Thus, $\md_1$ coincides with the de Rham operator on $M_1$.  Here
$\cP\in \Gamma (\Lambda^2 TM_2)$ and $\de\in \Gamma (TM_2)\otimes
\Omega^1 (M_2)$ are sections of $\Omega^\cdot_E (M)\otimes
\Omega^\cdot (M)$ living only over $M_2$.  Since $\int_{M_1} L_\xi$
equals zero for any vector field $\xi\in\Gamma(TM_1)$ (taking into
account that $L_\xi (\cdot )$ is always exact when acting on a form of
highest degree), it is sufficient to check the statement for an
arbitrary ``vertical'' section $\e\in
\Gamma (\mathrm{pr}^*_2 E_2)$ (whose projection to $TM_1$ vanishes).
One can easily calculate that
\beq
\EL_\e (\de) &=& \md_1 \e_i \otimes \md X^i + \cP^{ji} \pr_i \otimes \md_1
\e_j + (\e_{j,i}-\e_{i,j})\cP^{jk} \pr_k \otimes \md X^i
\in \Omega^1_E (M)\otimes \Omega^1 (M) \\
\EL_\e (\cP ) &=&  \cP^{ji} \md_1\e_j \otimes \pr_i +\frac{1}{2}
(\e_{j,i}-\e_{i,j})\cP^{ki}\cP^{lj} \pr_k\wedge \pr_l \in \Omega^2_E (M)
\eeq
which implies that the corresponding variation of the PSM action in
the form (\ref{actiondelta}) equals
\beq\label{PSMvar}
\EL_\e S_{PSM} =\int\limits_{M_1}
\md_1 \e_i \wedge \md X^i + (\e_{j,i}-\e_{i,j})\cP^{jk} A_k \otimes \md X^i +
\frac{1}{2} (\e_{j,i}-\e_{i,j})\cP^{ki}\cP^{lj} A_k\wedge A_l \, . 
\eeq
Clearly, the expression (\ref{PSMvar}) vanishes if $\pr M_1 =\emptyset$ and
 the required condition (\ref{gauge_condition}) holds, which
implies that $(\e_{j,i}-\e_{i,j})\cP^{jk}\pr_k \equiv 0$.

$\blacksquare$

\vskip7mm

In the remainder we briefly compare with another point of view on
gauge transformations, viewed as an action of a certain infinite
dimensional Lie algebroid living on the space of base maps,
c.f. \cite{LO}. Let $E_i$ be Lie algebroids over $M_i$, $i=1,2$. Then
there is a vector bundle $\cE$ over the space $\cM_0$ of smooth maps
$\phi_0$ acting from $M_1$ to $M_2$, defined such that the infinite
dimensional fiber $\cE_{\phi_0 }$ at any point $\phi_0$ is $\Gamma
(M_1 ,\phi_0^* E_2)$.

One has a natural map $\mathrm{Ind}^E$ acting from sections of
$\mathrm{pr}_2^* E_2$ over $M$, as used before, to sections of $\cE$
over $\cM$: any section $s\in \Gamma (M ,\mathrm{pr}_2^* E_2)$ gives a
section of $\cE$ by the map $s\mto \mathrm{Ind}^E_s$, such that $
\mathrm{Ind}^E_s (\phi_0) := (\phigra_0)^* s\in \Gamma (M_1 ,
\phi_0^* E_2)$. The map $\mathrm{Ind}^E$ is an embedding; moreover,
the space of all sections $\Gamma (\cM_0 ,\cE )$ is generated by
$\mathrm{Ind}^E_s$, $s\in \Gamma (M ,\mathrm{pr}_2^* E_2)$ over an
appropriate algebra of ``smooth'' functions on $\cM_0$.  For example,
if $E_2 =TM_2$ then the corresponding bundle over $\cM_0$ can be
thought of as $T\cM_0$. Let us notice that $T\cM_0$ is also a Lie
algebroid, such that the map $\mathrm{Ind}^T\colon
\Gamma (M, \mathrm{pr}_2^* TM_2 )\to \Gamma (\cM_0 ,T\cM_0)$
respects the Lie brackets.
We can easily extend this fact for a general $\cE\to \cM_0$ obtained as
above. For this purpose we introduce an anchor map $
{\,\!{}^{\cE}\!\rho} \colon \cE \to T\cM_0$ such that the following
diagram is commutative:
\beq\be{ccc}
\Gamma (M, \mathrm{pr}_2^* E_2) &
\stackrel{\mathrm{Ind}^E}{\longrightarrow} & \Gamma (\cM_0 , \cE)
\\
\rho \; \downarrow & &  {\,\!{}^{\cE}\!\rho} \; \downarrow \\
\Gamma (M, \mathrm{pr}_2^* TM_2)  &
\stackrel{\mathrm{Ind}^T}{\longrightarrow} & \Gamma (\cM_0 ,T\cM_0)
\ee\eeq
Now the Lie bracket on the image of $\Gamma (M, \mathrm{pr}_2^* E_2)$
can be extended to the space of all sections of $\Gamma (\cM_0 ,\cE)$ by
which it becomes a Lie algebroid bracket.

As an example, consider $M_2=pt$, $E_2$ a Lie algebra $\g$ with
trivial anchor map $\rho\equiv 0$. Then $\cM_0$ consists of only one
element, and $\cE=C^\infty (M_1 ,\g)$ is an infinite-dimensional Lie
algebra of ``multiloops''.

In the language of Poisson Sigma Models, or more generally in the
setting of Theorem \ref{theorem1}, $\delta_\e$ defines a gauge
transformation for any section $\e\in\Gamma(\cM_0,\cE)$.  The previous
discussion, however, only led to an action of $\Gamma(\cE)$ on base
maps $\phi_0\colon M_1\to M_2$ via the vector field
${\,\!{}^{\cE}\!\rho}(\epsilon)$ on $\cM_0$. More generally, all
vector fields $v_1\in\Gamma(TM_1)$ and $v_2\in\Gamma(TM_2)$ define
sections $\bar{v}_1$ and $\bar{v}_2$ of $\Gamma(\cM_0)$ which in a point
$\phi_0\in\cM_0$ take the value
\begin{eqnarray}
 \bar{v}_1(\phi_0)(x) &:=& \md\phi_0\circ v_1(x)\quad\mbox{ and} \label{vf}\\
 \bar{v}_2(\phi_0)(x) &:=& v_2\circ\phi_0(x)\,,
\end{eqnarray} 
respectively. Here, we use $T_{\phi_0}\cM_0\cong \Gamma(M_1,\phi_0^*
TM_2)$ such that a vector field on $\cM_0$ is defined by giving its
value $v(\phi_0)(x)\in\phi_0^*TM_2$ in a map $\phi_0$ and a point
$x\in M_1$. Both vector fields can be seen to generate left and right
compositions of diffeomorphisms on $M_1$ and $M_2$, respectively, with
maps in $\cM_0$. As such, those vector fields always commute with each
other. Sections of $E_1$ and $E_2$ then define vector fields on
$\cM_0$ through $\rho_1(\epsilon_1)\in\Gamma(TM_1)$ and
$\rho_2(\epsilon_2)\in\Gamma(TM_2)$.

This construction is clearly not general enough for our purposes. For
gauge transformations we need vector fields which act on the set of
bundle maps $E_1\to E_2$ (i.e. ``classical fields'') denoted as
$\cM$.  
This space $\cM$ is a bundle over $\cM_0$ with fiber over a point
$\phi_0\in \cM_0$ equal to $\Omega^1_{E_1}(M_1, \phi_0^* E_2)$.

Vector fields on $\cM$ suitable for gauge transformations can
advantageously be defined in the framework of infinite-dimensional
supergeometry (however, an advantage of our independent construction
is that we avoid infinite-dimensional supercomplications).  A vector
bundle $E\to M$ can be thought of as a $\Z-$graded manifold, denoted
as $E[1]$, with the parity of the fibers defined to be odd.  The
algebra of smooth functions $C^{\infty} (E [1])$ on $E[1]$ is
naturally isomorphic to $\Gamma (M ,\Lambda^\cdot E^*)$, and any
bundle map $E_1\to E_2$ between two vector bundles becomes a degree
preserving map $E_1 [1]\to E_2 [1]$. For any Lie algebroid $E\to M$
the canonical differential $\Emd$ defines a (super-) vector field of
degree one tangent to $E[1]$, endowing $E[1]$ with a $Q$-structure.
(A $\Z$-graded manifold endowed with an odd nilpotent vector field is
called a $Q$-manifold \cite{AKSZ}.)  Using this formalism, we can
reformulate the chain property (\ref{LA_morph}): a Lie algebroid
morphism is a map $\phi\colon E_1[1]\to E_2 [1]$ of degree zero, such
that $\phi_* (\md_1 )=\md_2$.

Denote the space of all graded maps $E_1 [1]\to E_2 [1]$ as $\cM_{\Z}$
(containing $\cM$ as the zero degree part).  Analogously to the
previous construction (\ref{vf}), the vector fields $\md_1$ and
$\md_2$ on $E_1 [1]$ and $E_2 [1]$ naturally generate commuting vector
fields $\de_1$ and $\de_2$ on $\cM_Z$, respectively (corresponding to
left and right compositions of morphisms). Since $\md_1$ and $\md_2$
are odd and nilpotent, so are $\de_1$ and $\de_2$. The difference $\de
:=\de_1 - \de_2$ is again a nilpotent vector field of degree
$1$. Moreover, $\de$ vanishes on the set of maps which preserve the
$Q$-structures (in particular, on the set of Lie algebroid morphisms).

A Lie algebroid $E$ can be identified with the tangent bundle $TE[1]$
where the action of a vector field on functions $C^{\infty}
(E[1])\cong\Gamma (M ,\Lambda^\cdot E^*)$ is obtained by contraction
between $E$ and $E^*$. If we have a section $\e\in \Gamma (\cM_0,
\cE)$ taking values in $E_2$, we obtain a vector field
$\bar{\epsilon}$ on $\cM$. Using $T_{\phi}\cM\cong
\Gamma(E_1,\phi^*TE_2)$, the vector field
$\bar{\epsilon}\in\Gamma(T\cM)$ is defined by
\[
 \bar{\epsilon}(\phi)(x):=\epsilon_{\phi_0}(\pi_1(x))\circ\phi(x)
\]
for $x\in E_1$ and with $\pi_1\colon E_1\to M_1$. Using the super
structure of $\cM_{\Z}$, $\bar{\epsilon}$ is a vector field of degree
$-1$.  A straightforward computation shows that the supercommutator
between $\de$ and the contraction with $\e$ is a degree preserving
vector field (therefore it is tangent to the subspace $\cM$). This
formula for a generalized gauge flow expressed as a supercommutator is
an analog of Cartan's magic formula (\ref{Lie_derivative}), which now
holds in the context of an infinite-dimensional geometry of graded
maps. One can use this infinitesimal transformation to generalize the
gauge transformation (\ref{newsymm}) to sections $\e$ which not only
depend on $X\in M_2$, but also depend on the map $\phi_0$
nontrivially. In particular, $\e$ might be a local functional
determined by higher jets of a base map $M_1\to M_2$. In a similar
way, we can express sections of $\epsilon_1\in\Gamma(E_1)$ as vector
fields on $\cM$:
\[
 \bar{\epsilon}_1(\phi)(x):= \phi\circ\epsilon_1(\pi_1(x))\,.
\]
Note that, unlike the vector fields defined in (\ref{vf}), vector
fields obtained in this way from $\epsilon_1\in\Gamma(E_1)$ and
$\epsilon\in\Gamma(\cM_0,\cE)$ do not commute in general since
$\epsilon$ also depends on $M_1$.


\section*{Acknowledgement}
T.S.~gratefully acknowledges helpful discussions with A.~Alekseev and
D.~Roytenberg in an early stage of this work.


\end{document}